\documentclass[12pt]{amsart}

\usepackage{amsmath}
\usepackage{amsfonts}
\usepackage{amssymb}
\usepackage{mathtools}
\usepackage{amscd}
\usepackage{graphicx}
\usepackage[dvipsnames]{xcolor}

\usepackage{booktabs}
\usepackage{multirow}\usepackage{subcaption}
\usepackage{here}

\newtagform{tiny}{\tiny(}{)}

\usepackage{hyperref}
\hypersetup{pdftitle={Classification of strongly F-regular and F-pure surface singularities with imperfect residue fields},pdfauthor={Tatsuro Kawakami and Kenta Sato}}
\usepackage[margin=1.25in]{geometry}

\usepackage{amsthm}
\usepackage{comment}
\usepackage[all,cmtip]{xy}
\usepackage{tikz-cd}
\usetikzlibrary{cd}

\usepackage{etoolbox}
\apptocmd{\sloppy}{\hbadness 10000\relax}{}{}

\usepackage{enumitem}

\title{Classification of strongly $F$-regular and $F$-pure surface singularities with imperfect residue fields}
\author{Tatsuro Kawakami}
\address{Graduate School of Mathematical Sciences, University of Tokyo, 3-8-1 Komaba,
Meguro-ku, Tokyo 153-8914, Japan}
\email{kawakami@ms.u-tokyo.ac.jp}

\author{Kenta Sato}
\address{Department of Mathematics and Informatics, Chiba University, Chiba, 263-8522, Japan} 
\email{sato@math.s.chiba-u.ac.jp}

\def\ge{\geqslant}
\def\le{\leqslant}
\def\phi{\varphi}
\def\epsilon{\varepsilon}
\def\tilde{\widetilde}

\def\mapsto{\longmapsto}

\def\Hom{\operatorname{Hom}}
\def\Spec{\operatorname{Spec}}
\def\Proj{\operatorname{Proj}}
\def\Supp{\operatorname{Supp}}

\def\Exc{\operatorname{Exc}}

\def\max{\operatorname{max}}

\def\id{\mathrm{id}}

\def\m{{\mathfrak m}}

\newcommand{\Q}{\mathbb{Q}} 
\newcommand{\C}{\mathbb{C}} 
\newcommand{\R}{\mathbb{R}} 
\newcommand{\Z}{\mathbb{Z}}

\newcommand{\PP}{\mathbb{P}}

\newcommand{\sO}{\mathcal{O}}

\newcommand{\stacksproj}[1]{{\cite[\href{https://stacks.math.columbia.edu/tag/#1}{Tag~{#1}}]{stacks-project}}}

\newsavebox{\pullback}
\sbox\pullback{%
\begin{tikzpicture}%
\draw (0,0) -- (1ex,0ex);%
\draw (1ex,0ex) -- (1ex,1ex);%
\end{tikzpicture}}

\newsavebox{\pullbackdl}
\sbox\pullbackdl{%
\begin{tikzpicture}%
\draw (-1ex,0ex) -- (0ex,0ex);%
\draw (0ex,-1ex) -- (0ex,0ex);%
\end{tikzpicture}}

\newsavebox{\pushoutdr}
\sbox\pushoutdr{%
\begin{tikzpicture}%
\draw (-1ex,-1ex) -- (-1ex,0ex);%
\draw (-1ex,0ex) -- (0ex,0ex);%
\end{tikzpicture}}

\theoremstyle{plain}
\newtheorem{thm}{Theorem}[section] 
\newtheorem{cor}[thm]{Corollary}
\newtheorem{prop}[thm]{Proposition}

\newtheorem{lem}[thm]{Lemma}
\theoremstyle{definition} 
\newtheorem{defn}[thm]{Definition}

\newtheorem{setting}[thm]{Setting}
 
\newtheorem{eg}[thm]{Example} 

\theoremstyle{remark}
\newtheorem{rem}[thm]{Remark}

\newtheorem{defn and notation}[thm]{Definition and Notation}

\newtheorem*{cl}{Claim}

\theoremstyle{plain}
\newtheorem{theo}{Theorem}

\numberwithin{equation}{thm}
\keywords{$F$-singularities; surface singularities; imperfect fields; dual graphs}
\subjclass[2020]{13A35, 14B05, 14J17}

\begin{document}
\tolerance = 9999

\begin{abstract}
We study strongly $F$-regular and $F$-pure normal surface singularities whose residue fields are not necessarily perfect. 
We characterize strong $F$-regularity in terms of the weighted dual graph and geometric reducedness of the exceptional components of the minimal resolution, extending Hara's classification over algebraically closed residue fields. 
We also determine the possible dual graphs of $F$-pure surface singularities and give converse criteria. 
Notably, the exceptional components of an $F$-pure surface singularity are geometrically reduced unless the singularity is a rational double point. 
\end{abstract}

\maketitle
\markboth{Tatsuro Kawakami and Kenta Sato}{STRONGLY F-REGULAR AND F-PURE SURFACE SINGULARITIES}

\section{Introduction}\label{sec:introduction}

Strongly $F$-regular and $F$-pure singularities are fundamental classes of singularities in positive characteristic. For normal $\Q$-Gorenstein singularities, they imply Kawamata log terminal (klt) and log canonical (lc) singularities, respectively \cite{Hara-Watanabe}. In dimension two, Hara \cite{Hara(two-dim)} classified strongly $F$-regular and $F$-pure singularities with algebraically closed residue fields in terms of their minimal resolutions. In particular, a klt surface singularity with algebraically closed residue field is strongly $F$-regular in characteristic $p>5$. In smaller characteristics, the dual graph imposes additional restrictions on $p$.

The assumption on the residue field is a substantial restriction in applications to higher-dimensional varieties. For example, at a codimension-two point of a threefold over an algebraically closed field, the local ring is two-dimensional but its residue field has transcendence degree one and is imperfect. Thus, even over an algebraically closed base field, arguments involving localization in codimension two require results over imperfect residue fields. Sato and Takagi \cite{Sato-Takagi(generalhyperplane)} proved that a two-dimensional Noetherian $F$-finite normal local domain of characteristic $p>5$ is klt if and only if it is strongly $F$-regular, without assuming that its residue field is perfect. Their result also yields a Bertini theorem for klt threefolds in this range of characteristics.

In this paper, we extend Hara's classification to arbitrary residue fields, including characteristics $2$, $3$, and $5$. Two features distinguish this setting from the algebraically closed case. Let $\pi\colon T\to \Spec R$ be the minimal resolution of a klt surface singularity with residue field $\kappa$, and write $\Exc(\pi)=\bigcup_i E_i$. Each $E_i$ has genus zero, but it need not be smooth, or even geometrically reduced, over $\kappa$. Moreover, the fields $H^0(E_i,\sO_{E_i})$, and the residue fields of intersection points, need not equal $\kappa$. The dual graph must therefore record the corresponding field degrees as well as genera and intersection numbers. We use the weighted graphs and the terminology of twisted quotients in \cite[Section 3.3]{Kol13} and \cite[Appendix A]{Sato23} (see Section~\ref{sec:preliminaries}).

Our first main result shows that geometric reducedness of the exceptional components, together with the same restrictions on graph types and characteristic as in Hara's classification, characterizes strong $F$-regularity.

\begin{theo}\label{Introthm:classification of SFR}
Let $(R,\m,\kappa)$ be a two-dimensional Noetherian $F$-finite normal klt local domain of characteristic $p>0$, and let $\pi\colon T\to S\coloneqq\Spec R$ be its minimal resolution, with exceptional components $E_1, \dots, E_n$.
Then $R$ is strongly $F$-regular if and only if either $R$ is regular or both of the following conditions hold:
\begin{enumerate}[label=\textup{(\Alph*)}]
\item Every $E_i$ is geometrically reduced over $\kappa$.
\item The dual graph of $\pi$ is one of the following:
\begin{enumerate}[label=\textup{(\arabic*)}]
\item A chain or a twisted chain.
\item A star-shaped or twisted star-shaped graph of one of the following types:
\begin{enumerate}[label=\textup{(\roman*)}]
\item $(2,2,d)$ with $d\geq2$ and $p>2$;
\item $(2,3,3)$ or $(2,3,4)$ with $p>3$;
\item $(2,3,5)$ with $p>5$.
\end{enumerate}
\end{enumerate}
\end{enumerate}
\end{theo}

If $p>3$, condition~(B) already implies condition~(A) (Corollary~\ref{cor:geom red,p>3}). In particular, for $p>5$, Theorem~\ref{Introthm:classification of SFR} recovers the equivalence between klt and strongly $F$-regular surface singularities proved in \cite{Sato-Takagi(generalhyperplane)}. In characteristics $2$ and $3$, geometric reducedness is an additional arithmetic condition that is not determined by the weighted dual graph alone.
See Example \ref{eg:geom reducedness of exc div2} for more details.

Our second main result concerns $F$-pure singularities. 
Since $F$-pure singularities are lc, the possible graphs are the following:
the star-shaped graphs of type $(a,b,c)$ with $1/a+1/b+1/c \ge 2$ , the graphs $\tilde D_n$ with $n \ge 4$, cusps, simple elliptic singularities, and their twisted forms.
We use ``simple elliptic'' in the sense of \cite[Definition 4.12]{Sato23}, which includes the simple nodal elliptic case.

\begin{theo}\label{Introthm:classification of SFP}
Let $(R,\m,\kappa)$ be a two-dimensional Noetherian $F$-finite normal lc local domain of characteristic $p>0$.
\begin{enumerate}
\item[\textup{(I)}] If $R$ is $F$-pure, then at least one of the following holds:
\begin{enumerate}[label=\textup{(\arabic*)}]
\item $R$ is strongly $F$-regular.
\item $R$ is a rational double point.
\item The dual graph is star-shaped or twisted star-shaped of one of the following types:
\begin{enumerate}[label=\textup{(\roman*)}]
\item $(3,3,3)$ or $(2,3,6)$, with $p\equiv1\pmod3$;
\item $(2,4,4)$, with $p\equiv1\pmod4$;
\item $(2,2,2,2)$ (equivalently, $\tilde D_4$), with $p>2$.
\end{enumerate}
\item The dual graph is $\tilde D_n$ or twisted $\tilde D_n$ for some $n\geq5$, and $p>2$.
\item $R$ is a cusp or a twisted cusp.
\item $R$ is a simple elliptic singularity.
\end{enumerate}
\item[\textup{(II)}] Conversely, each of \textup{(1)}, \textup{(3-i)}, \textup{(3-ii)}, and \textup{(4)} implies that $R$ is $F$-pure. Condition~\textup{(5)} also implies $F$-purity if $p>2$.
\end{enumerate}
\end{theo}

The cases omitted from the converse require additional information. For rational double points, both $F$-pure and non-$F$-pure examples occur, as over an algebraically closed residue field. Such singularities are hypersurfaces \cite[Proposition 3.9]{Sato-Takagi(generalhyperplane)}, so their $F$-purity can be tested by Fedder's criterion \cite{Fed}. For a simple elliptic singularity, $F$-purity is equivalent to global $F$-splitting of the exceptional curve. For the $\tilde D_4$ configurations, the corresponding criterion is global $F$-splitting of a genus-zero curve with its induced boundary; in the non-twisted case this gives an explicit coefficient test (Remark~\ref{rem:Remark on Theorem B} and Example~\ref{eg:Theorem B for D4}).

Twisted cusps exhibit a further phenomenon in characteristic $2$. Although cusp singularities with algebraically closed residue field are $F$-pure in every characteristic, twisted cusps over imperfect fields need not be $F$-pure. Example~\ref{eg:non F-pure cusp} gives an explicit example. More precisely, a cusp or twisted cusp is $F$-pure if and only if all exceptional components of its minimal resolution are geometrically reduced; the extra condition is automatic when $p>2$ (Remark~\ref{rem:Remark on Theorem B}).

The main geometric input in both classifications is the following statement.

\begin{theo}[Corollary~\ref{cor:SFP to Geom red}]
Let $(R,\m,\kappa)$ be a two-dimensional Noetherian $F$-finite normal local domain of characteristic $p>0$. If $R$ is $F$-pure and $\Spec R$ is not a rational double point, then every exceptional component of its minimal resolution is geometrically reduced over $\kappa$.
\end{theo}

For strongly $F$-regular singularities, the same conclusion holds without the exception for rational double points (Theorem~\ref{thm:SFR to Geom.red}). To prove these assertions, we relate geometric reducedness to separability of the constant fields of exceptional components and use global $F$-splitting on suitable exceptional curves. 

Once geometric reducedness has been proven, 
we can prove Theorem A and Theorem B by an argument similar to that in \cite{Hara(two-dim)}, except in the case of rational double points.
In fact, after finite separable extensions of the residue field, the genus-zero configurations in the dual graphs become configurations of projective lines, to which Hara's method applies.

Rational double points require a separate argument. Over an algebraically closed field, Hara uses Artin's equations \cite{Art} and Fedder's criterion to obtain the restrictions on the characteristic. 
We instead use an intrinsic consequence of strong $F$-regularity in the Gorenstein case: there is a nonzero $f\in\m$ such that $(\Spec R,\frac{1}{p-1}\operatorname{div}(f))$ is sharply $F$-pure (Lemma~\ref{lem:key}). 
The coefficients of the pullback of $\operatorname{div}(f)$ satisfy intersection inequalities on the minimal resolution. 
These inequalities yield the required lower bounds on $p$ without explicit equations for the singularity.

The classification has important applications to differential forms in higher dimensions. We can apply this classification to prove the logarithmic extension theorem for one-forms on strongly $F$-regular singularities in arbitrary dimension established in \cite{Kawakami-Sato}.

The paper is organized as follows. Section~\ref{sec:preliminaries} recalls the required facts about $F$-singularities and dual graphs of curves on surfaces. Section~\ref{sec:geometric-reducedness} studies geometric reducedness of exceptional curves. Section~\ref{section:Appendix} proves the main theorems.

\section{Preliminaries}\label{sec:preliminaries}

\subsection{Notation and terminology}
Throughout this paper, we use the following notation:
\begin{itemize}
\item We say that an $\mathbb{F}_p$-scheme $X$ is {\em $F$-finite} if $F \colon X \to X$ is a finite morphism, and  we say that an $\mathbb{F}_p$-algebra $R$ is {\em $F$-finite} if $\Spec R$ is $F$-finite. If $R$ is a Noetherian $F$-finite $\mathbb{F}_p$-algebra, it is excellent, admits a dualizing complex, and $\dim\,R<\infty$ (cf.~\cite[Remark 13.6]{Gab04}).
\item Given a proper birational morphism $f\colon Y\to X$, we denote by $\Exc(f)$ the exceptional locus with reduced structure.
\item For the definition of the singularities of pairs appearing in the MMP (such as \emph{klt, plt, lc}), we refer to \cite[Definition 2.8]{Kol13}. Note that we always assume that the boundary divisor is effective although \cite[Definition 2.8]{Kol13} does not impose this assumption.
\item Given a property $\mathcal{P}$ of schemes, we say that a scheme $X$ over a field $k$ is \emph{geometrically $\mathcal{P}$} if for every finite field extension $K$ over $k$, the base change $X_K \coloneqq X \times_k \Spec K$ satisfies $\mathcal{P}$.
\end{itemize}

\subsection{Sharp \texorpdfstring{$F$}{F}-splitting and strong \texorpdfstring{$F$}{F}-regularity}

In this subsection, we recall basic properties of sharp $F$-splitting and strong $F$-regularity.

\begin{defn}\label{def:GFS and GFR}
Let $X$ be a normal $F$-finite scheme of characteristic $p>0$, and let $D$ be an effective $\Q$-divisor on $X$. 
\begin{enumerate}
\item[\textup{(1)}] We say that $(X, D)$ is {\em globally sharply $F$-split} if 
\[
\sO_X \xrightarrow{F^e} F_*^e\sO_X \hookrightarrow F_*^e\sO_X( \lceil(p^e-1)D \rceil)
\]
splits as an $\sO_X$-module homomorphism for some $e \in \Z_{>0}$. 
 
\item[\textup{(2)}] 
We say that $(X, D)$ is {\em globally $F$-regular} 
if, for every effective $\Z$-divisor $E$, there exists $e \in \Z_{>0}$ such that 
\[
\sO_X \xrightarrow{F^e} F_*^e\sO_X \hookrightarrow 
F_*^e\sO_X( \lceil(p^e-1)D \rceil +E)
\]
splits as an $\sO_X$-module homomorphism. 
\end{enumerate}
\end{defn}

For sharp $F$-splitting with $D=0$, we simply say that $X$ is \emph{globally $F$-split}.
More generally, this term is used for any $F$-finite scheme $X$, without a normality assumption, if $\sO_X\to F_*\sO_X$ splits.
In particular, it applies to reduced exceptional curves that need not be normal.

\begin{lem}\label{lem:sharply or not}
    Let $X$ be a normal variety over an $F$-finite field $k$ of characteristic $p>0$, and let $D$ be an effective $\Q$-divisor with standard coefficients, that is, coefficients in $\{\frac{m-1}{m}\mid m\in\Z_{>0}\}$.
    Then the pair $(X,D)$ is globally sharply $F$-split if and only if it is $F$-split in the sense of \cite[Definition 2.4]{Hara(two-dim)}, that is, the natural morphism
    \[
    \sO_X \to F_{*}\sO_X(\lfloor pD \rfloor)
    \]
    splits as $\sO_X$-module homomorphism.
\end{lem}

\begin{proof}
For every positive integer $q$ and every $m\geq 1$, one has
\[
\left\lceil(q-1)\frac{m-1}{m}\right\rceil
=\left\lfloor q\frac{m-1}{m}\right\rfloor.
\]
Thus a splitting for $e=1$ in the displayed map gives sharp $F$-splitting.
Conversely, suppose that the map for $q=p^e$ splits.
For a standard coefficient $a$, the inequality
$p^{e-1}\lfloor pa\rfloor\leq\lfloor p^ea\rfloor$ holds.
Consequently, the splitting factors through
\[
\sO_X\longrightarrow F_*\sO_X(\lfloor pD\rfloor)
\longrightarrow F_*^e\sO_X(\lfloor p^eD\rfloor),
\]
where the second arrow is induced by the $(e-1)$st Frobenius.
The first arrow therefore splits.
\end{proof}

\begin{defn}\label{def:GFS and GFR(local)}
Let $X$ be a normal $F$-finite scheme of characteristic $p>0$ and let $D$ be an effective $\Q$-divisor on $X$.
\begin{enumerate}
\item[\textup{(1)}] We say that the pair $(X, D)$ is {\em sharply $F$-pure} if, for every $x\in X$, the pair $(\Spec\,\sO_{X,x},D_x)$
is globally sharply $F$-split.
When we have $D=0$, we simply say that $X$ is \emph{$F$-pure}.
\item[\textup{(2)}] We say that the pair $(X, D)$ is {\em strongly $F$-regular} if, for every $x\in X$, the pair $(\Spec\,\sO_{X,x},D_x)$
is globally $F$-regular.
\end{enumerate}
\end{defn}

\begin{rem}
All the properties in Definition \ref{def:GFS and GFR(local)} can be checked at closed points.
When $X$ is affine, sharp $F$-purity (resp.~strong $F$-regularity) is equivalent to global sharp $F$-splitting (resp.~global $F$-regularity). 
\end{rem}

The following lemma is one of the key ingredients for the proof of Theorem \ref{Introthm:classification of SFR}.

\begin{lem}\label{lem:key}
    In Setting \ref{notation:normal surface sings}, suppose that $R$ is strongly $F$-regular and Gorenstein.
    Then there is a nonzero element $f\in\m$ such that
    $(X,\frac{1}{p-1}D)$ is sharply $F$-pure, where $D=\operatorname{div}(f)$.
\end{lem}
\begin{proof}
Since $R$ is Gorenstein, $\Hom_R(F_*R,R)$ is free of rank one as an $F_*R$-module. We fix a generator $\phi$.
If $\phi(F_*\m)\subseteq\m$, then $\m$ is compatible with every Frobenius map, contradicting strong $F$-regularity
(see \cite[Propositions 4.1 and 4.8(i)]{Schwede(F-adj)}).
Thus $\phi(F_*\m)=R$, and we can choose $f\in\m$ with $\phi(F_*f)=1$.
Since $R(D)=f^{-1}R$, multiplication by $f$ gives an isomorphism
$F_*R(D)\xrightarrow{\sim}F_*R$.
Its composite with $\phi$ sends $F_*1$ to $1$, and hence splits
$R\to F_*R(D)$.
\end{proof}

We will use adjunction for sharp $F$-splitting:

\begin{lem}\label{lem:adjunction}
    Let $X$ be a normal integral Noetherian $F$-finite scheme of characteristic $p>0$ with a dualizing complex.
    Let $S$ be a normal prime divisor and let $B$ be an effective $\Q$-divisor such that
    $S\not\subseteq\Supp B$ and $K_X+S+B$ is $\Q$-Cartier.
    Write
    \[
    (K_X+S+B)|_S=K_S+B_S,
    \]
    where $B_S$ is the different.
    If $(X,S+B)$ is globally sharply $F$-split, then so is $(S,B_S)$.
\end{lem}
\begin{proof}
Using global splitting of $(X,S+B)$, we can find an effective $\Q$-divisor
$S+B'\geq S+B$ such that $(X,S+B')$ is globally sharply $F$-split and
$(p^e-1)(K_X+S+B')$ is Cartier for some $e>0$.
(see \cite[proof of Theorem 4.3(ii)]{SS10}).
The coefficient of $S$ remains one, since sharp $F$-purity forces every boundary coefficient to be at most one.
Let $B'_S$ be the different of $B'$ on $S$.
By the compatibility of Frobenius trace, together with the equality of the different and the $F$-different
\cite[Theorem 5.3]{Das15}, we have a commutative diagram
\[
\begin{tikzcd}
H^0\bigl(X,F_*^e\sO_X((1-p^e)(K_X+S+B'))\bigr)
 \arrow[r,"\operatorname{Tr}"]\arrow[d]
 & H^0(X,\sO_X)\arrow[d]\\
H^0\bigl(S,F_*^e\sO_S((1-p^e)(K_S+B'_S))\bigr)
 \arrow[r,"\operatorname{Tr}"]
 & H^0(S,\sO_S).
\end{tikzcd}
\]
A section mapping to $1$ in the upper row restricts to a section mapping to $1$ in the lower row.
Thus $(S,B'_S)$ is globally sharply $F$-split.
Since $B'-B$ is effective and $\Q$-Cartier and does not contain $S$, adjunction gives
$B'_S-B_S=(B'-B)|_S\geq0$.
Sharp $F$-splitting is preserved when the boundary decreases, so $(S,B_S)$ is globally sharply $F$-split as well.
\end{proof}

\begin{rem}\label{rem:non-normal adjunction}
Let $X$ be as in Lemma \ref{lem:adjunction} and let $E$ be a (not necessarily normal) reduced Cartier divisor on $X$.
By an argument similar to that in Lemma \ref{lem:adjunction}, if $(X,E)$ is globally sharply $F$-split, then so is $E$.
\end{rem}

\subsubsection{\texorpdfstring{$F$}{F}-splitting and base changes}

\begin{thm}\label{thm:F-split to geom F-split}
    Let $K/k$ be a field extension of $F$-finite fields.
    Let $X$ be a geometrically normal proper variety over $k$, and let $D$ be an effective $\Q$-divisor on $X$. Assume that $H^0(X,\sO_X)$ is a separable field extension of $k$.
    If $(X,D)$ is globally sharply $F$-split (resp.~globally $F$-regular), then so is $(X\times_k K,D\times_k K)$.
\end{thm}
\begin{proof}
For global $F$-regularity, this is \cite[Corollary 2.8]{GLP$^+$15}.
For sharp $F$-splitting, choose $a>0$ for which the defining map splits and replace
$D$ by $\lceil(p^a-1)D\rceil/(p^a-1)$.
We can apply the proof of \cite[Proposition 2.6 and Remark 2.7]{GLP$^+$15}, with its auxiliary effective divisor set equal to zero, to this splitting. Here, we use geometric normality.
Then the argument in the proof of \cite[Corollary 2.8]{GLP$^+$15} gives the assertion for arbitrary $K/k$.
\end{proof}

\begin{lem}\label{lem:base change of GFR}
    Let $K/k$ be an extension of $F$-finite fields of characteristic $p>0$.
    Let $X$ be a normal variety over $k$, let $D$ be an effective $\Q$-divisor on $X$, and suppose that $X_K=X\times_kK$ is normal.
    Write $D_K$ for the flat pullback of $D$.
    \begin{enumerate}[label=\textup{(\arabic*)}]
    \item If $(X_K,D_K)$ is globally $F$-regular, then so is $(X,D)$.
    \item Suppose, in addition, that every prime component of $D$ has reduced pullback to $X_K$.
    If $(X_K,D_K)$ is globally sharply $F$-split, then so is $(X,D)$.
    \end{enumerate}
\end{lem}
\begin{proof}
Let $\rho\colon X_K\to X$ be the projection, and choose a $k$-linear map
$\lambda\colon K\to k$ with $\lambda(1)=1$.
The map $\mathrm{id}\otimes\lambda$ is an $\sO_X$-linear retraction
$\rho_*\sO_{X_K}\to\sO_X$.

For (1), fix an effective integral divisor $E$ on $X$ and set $F=D_{\mathrm{red}}$.
Global $F$-regularity provides $e>0$ and a splitting
\begin{equation}\label{eq:gfr}
\sO_{X_K}\longrightarrow
F_*^e\sO_{X_K}\bigl(\lceil(p^e-1)D_K\rceil+(E+F)_K\bigr).
\end{equation}
Since
\[
\bigl(\lceil(p^e-1)D\rceil\bigr)_K
\leq\lceil(p^e-1)D_K\rceil+F_K,
\]
there is a natural $\sO_X$-linear map
\[
F_*^e\sO_X\bigl(\lceil(p^e-1)D\rceil+E\bigr)
\longrightarrow
\rho_*F_*^e\sO_{X_K}\bigl(\lceil(p^e-1)D_K\rceil+(E+F)_K\bigr).
\]
Composing this map with the retraction in \eqref{eq:gfr} and then with
$\mathrm{id}\otimes\lambda$ gives a map to $\sO_X$ sending $1$ to $1$.
This proves (1).

For (2), the assumption gives
\[
\bigl(\lceil(p^e-1)D\rceil\bigr)_K
=\lceil(p^e-1)D_K\rceil.
\]
The same argument, omitting the auxiliary divisors $E$ and $F$, descends a sharp splitting.
\end{proof}

\begin{rem}
The extra hypothesis in part (2) is essential for purely inseparable extensions.
For example, let $k=\mathbb F_p(a)$, $K=k(a^{1/p})$, and
$P=V(t^p-a)\subset\PP^1_k$.
For $D=[\infty]+\frac1pP$, the pair $(\PP^1_K,D_K)$ is globally sharply $F$-split,
since $D_K$ consists of two distinct $K$-rational points with coefficient one.
On the other hand, for $q=p^e$ one has
$\deg\lceil(q-1)D\rceil=2q-1$.
Thus $\sO_{\PP^1_k}((1-q)K_{\PP^1_k}-\lceil(q-1)D\rceil)$ has degree $-1$,
so Frobenius duality rules out a global sharp splitting over $k$.
\end{rem}

The following classification \cite{Watanabe(one-dim)} remains valid over an arbitrary $F$-finite field when the marked points are rational.

\begin{prop}\label{prop:log Fano curve}
    Let $D$ be a $\Q$-divisor on $\mathbb{P}_k^1$ over an $F$-finite field $k$ of characteristic $p>0$ of the form
    \[
    D=\sum_{i=1}^r \frac{d_i-1}{d_i}P_i,
    \]
    where the $P_i$ are distinct $k$-rational points and $2 \leq d_1 \leq d_2 \leq \dots \leq d_r$.
    We denote $D$ by $(P_1, \dots, P_r; d_1,\ldots,d_r)$ and, when $r\leq3$, abbreviate this to $(d_1,\ldots,d_r)$, since the positions of the points do not affect the isomorphism class of the pair.
    If $r\leq2$, the pair $(\mathbb P^1_k,D)$ is globally $F$-regular in every characteristic.
    If $r\geq3$, it is globally $F$-regular if and only if one of the following holds:
    \begin{enumerate}
        \item[\textup{(1)}] $D=(2,2,d)$, $d\geq 2$, and $p\neq 2$.
        \item[\textup{(2)}] $D=(2,3,3)$ or $(2,3,4)$ and $p>3$.
        \item[\textup{(3)}] $D=(2,3,5)$ and $p>5$.
    \end{enumerate}
    Moreover, $(\mathbb{P}^1_k,D)$ is globally sharply $F$-split if and only if it is globally $F$-regular or one of the following holds:
    \begin{enumerate}
        \item[\textup{(4)}] $D=(3,3,3)$ or $(2,3,6)$, and $p\equiv 1$ modulo $3$.
        \item[\textup{(5)}] $D=(2,4,4)$ and $p\equiv 1$ modulo $4$.
        \item[\textup{(6)}] $D=(\infty,0,-1,\lambda;2,2,2,2)$ with $\lambda\in k\setminus\{-1,0\}$, and $p=2n+1$ such that the coefficient of $x^n$ in $(x+1)^n(x-\lambda)^n$ is nonzero.
    \end{enumerate}
\end{prop}
\begin{proof}
    Since the marked points are $k$-rational, their pullbacks remain reduced. Theorem \ref{thm:F-split to geom F-split} and Lemma \ref{lem:base change of GFR} therefore reduce the assertion to the case where $k$ is algebraically closed.
    For $r\leq2$, global $F$-regularity follows from the standard toric splitting of $\PP^1$ (or by a direct monomial computation). For $r\geq3$, the assertion follows from Lemma \ref{lem:sharply or not} and \cite[Proposition 2.6]{Hara(two-dim)} (cf.~\cite{Watanabe(one-dim)}).

\end{proof}

\subsection{Surface singularities}

\begin{setting}\label{notation:normal surface sings}
Let $(R,\m,\kappa)$ be a two-dimensional Noetherian $F$-finite normal local domain of characteristic $p>0$, and let $\pi\colon Y\to X \coloneqq \Spec R$ be the minimal resolution of $X$, with $\Exc(\pi)=\sum_{i=1}^n E_i$.
We further assume that $X$ has only lc singularities.
We denote by $K_i \coloneqq H^0(E_i,\sO_{E_i})$ and by $r_i$ the dimension $\dim_{\kappa} K_i$.
\end{setting}

In this paper, a graph is a $\mathbb{Z}^3$-weighted undirected multigraph.
With the above notation, the \emph{dual graph} of $X=\Spec R$ is a graph whose set of vertices is $\{E_1, \dots, E_n\}$, where the number of edges between $E_i$ and $E_j$ is the intersection number $E_i \cdot E_j$ ($= \dim_{\kappa} \mathcal{O}_{E_i \cap E_j}$), and the weight at $E_i$ is $(r_i, b_i \coloneqq E_i^2, g_i \coloneqq g(E_i)) \in \mathbb{Z}^3$, where $g(E_i) = \dim_{\kappa} H^1(E_i, \mathcal{O}_{E_i}) / r_i$ is the arithmetic genus of $E_i$.

Since $X$ is log canonical, every coefficient of $K_{Y/X}$ is larger than or equal to $-1$, that is, $X$ is \emph{numerically log canonical in the sense of \cite[Paragraph 2.27]{Kol13}}.
As is remarked in \cite[Paragraph 2.27]{Kol13}, the definition of numerically lc here is slightly different than the usual terminology.
The dual graphs of numerically lc surface singularities in this sense are classified in \cite[Subsection 3.3]{Kol13} (see also \cite[Appendix A]{Sato23}).
From this classification table,
we use the terms \emph{(twisted) chain}, \emph{(twisted) star-shaped graph}, \emph{(twisted) $\tilde D_n$}, \emph{(twisted) cusp} and \emph{simple elliptic}.
This terminology is numerical: a twisted graph need not be the graph obtained from the geometric exceptional divisor by Galois descent.
The following examples illustrate this distinction.

\begin{eg}
\
    \begin{enumerate}[label=\textup{(\arabic*)}]
    \item The minimal resolution of $X= \Spec \R[x,y,z]/(x^2+y^2+z^4)$ is obtained by blowing up twice.
    There are two exceptional curves $E_1$ and $E_2$, which satisfy 
    \[
    E_1 \cong \PP^1_{\C}, \ E_2 \cong \PP^1_{\R},\  (E_1^2)=-4,\ (E_2^2)=-2,\  \textup{and}\  E_1 \cap E_2 \cong \Spec \C.
    \]
    Therefore, the dual graph is as follows:
       
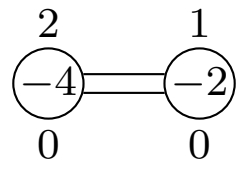
\[
\begin{tikzpicture}
\node[draw, shape=circle, inner sep=0.3pt] (E1) at (3,0){\tiny $-4$};
\node[draw, shape=circle, inner sep=0.3pt] (E2) at (4,0){\tiny $-2$};
\draw[shift={(0,0.4)}] (3,0) node{\tiny $2$};
\draw[shift={(0,0.4)}] (4,0) node{\tiny $1$};
\draw[shift={(0,-0.4)}] (3,0) node{\tiny $0$};
\draw[shift={(0,-0.4)}] (4,0) node{\tiny $0$};

\draw (E1.15)--(E2.165);
\draw (E1.-15)--(E2.-165);
\end{tikzpicture}\]   
This is the twisted chain associated with the following non-twisted graph:
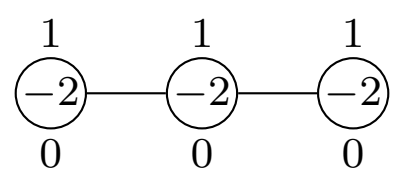
\[
\begin{tikzpicture}
\node[draw, shape=circle, inner sep=0.3pt] (E1) at (3,0){\tiny $-2$};
\node[draw, shape=circle, inner sep=0.3pt] (E2) at (4,0){\tiny $-2$};
\node[draw, shape=circle, inner sep=0.3pt] (E3) at (5,0){\tiny $-2$};

\draw[shift={(0,0.4)}] (3,0) node{\tiny $1$};
\draw[shift={(0,0.4)}] (4,0) node{\tiny $1$};
\draw[shift={(0,0.4)}] (5,0) node{\tiny $1$};
\draw[shift={(0,-0.4)}] (3,0) node{\tiny $0$};
\draw[shift={(0,-0.4)}] (4,0) node{\tiny $0$};
\draw[shift={(0,-0.4)}] (5,0) node{\tiny $0$};

\draw (E1)--(E2);
\draw (E2)--(E3);
\end{tikzpicture}\]   

We note that this graph is the dual graph of the \'etale cover
\[
X \times_\R \Spec \C=\Spec \C[x,y,z]/(x^2+y^2+z^4)
\]
of $X$.

\item We next consider $X' = \Spec \R[x,y,z]/(x^2+y^2+z^5)$, whose minimal resolution is obtained by blowing up twice.
    There are two exceptional curves $E_1$ and $E_2$, which satisfy 
    \[
    E_1 \cong \PP^1_{\C},  (E_1^2)=-4,\ (E_2^2)=-2,\  \textup{and}\  E_1 \cap E_2 \cong \Spec \C.
    \]
    Moreover, $E_2$ is isomorphic to the singular conic curve $Z(x^2+y^2) \subseteq \PP^2_{\R}$.
    Therefore, the dual graph is as follows:
       
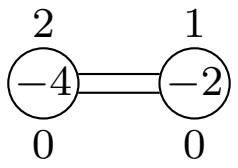
\[
\begin{tikzpicture}
\node[draw, shape=circle, inner sep=0.3pt] (E1) at (3,0){\tiny $-4$};
\node[draw, shape=circle, inner sep=0.3pt] (E2) at (4,0){\tiny $-2$};
\draw[shift={(0,0.4)}] (3,0) node{\tiny $2$};
\draw[shift={(0,0.4)}] (4,0) node{\tiny $1$};
\draw[shift={(0,-0.4)}] (3,0) node{\tiny $0$};
\draw[shift={(0,-0.4)}] (4,0) node{\tiny $0$};

\draw (E1.15)--(E2.165);
\draw (E1.-15)--(E2.-165);
\end{tikzpicture}
\]
In this case, the dual graph of $X'$ cannot be expressed as a twisted quotient of the dual graph of the \'etale cover 
\[
X' \times_{\R} \Spec \C = \Spec \C[x,y,z]/(x^2+y^2+z^5).
\]

\end{enumerate}    
\end{eg}

We use the notation of Setting \ref{notation:normal surface sings}.
If $r_i = 1$ for every $i$, we say that the dual graph of $X$ is \emph{non-twisted}.

\begin{rem}
Following \cite[Appendix A]{Sato23}, we use the term \lq\lq cusp\rq\rq for the following graphs, even though they are not necessarily non-twisted.
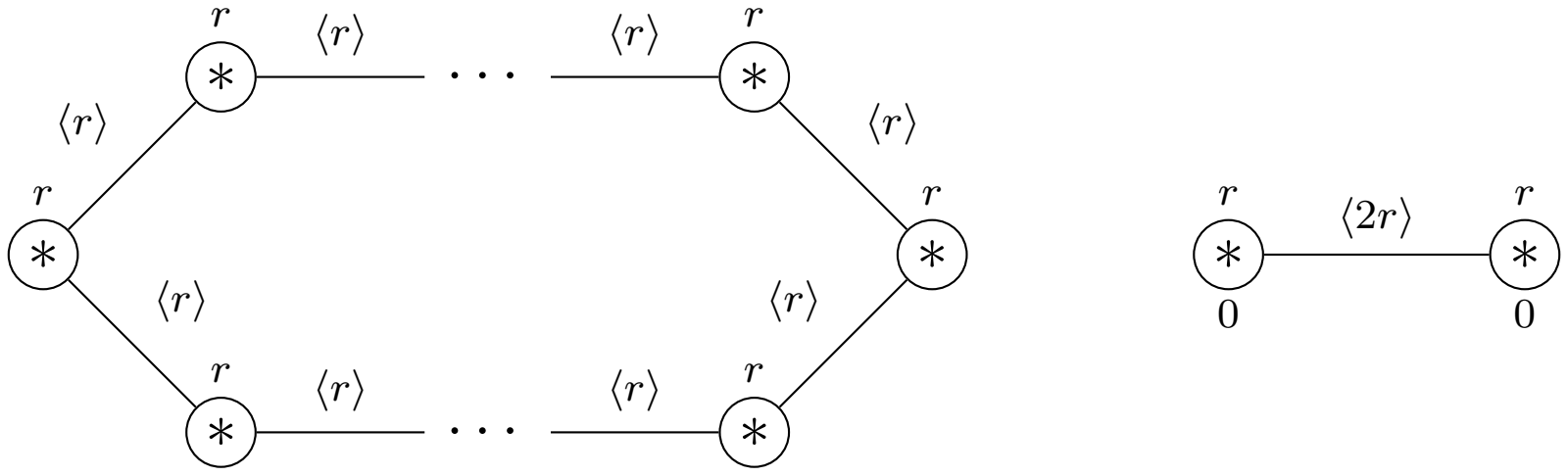
\[
\begin{tikzpicture}
\node[draw, shape=circle, inner sep=1.8pt] (A) at (0,0){$*$};
\draw[shift={(0,0.4)}] (0,0) node{\tiny $r$};
\node[draw, shape=circle, inner sep=1.8pt] (B) at (1.2,1.2){$*$};
\draw[shift={(0,0.4)}] (1.2,1.2) node{\tiny $r$};
\draw node (C) at (3,1.2){$\cdots$};
\node[draw, shape=circle, inner sep=1.8pt] (D) at (4.8 ,1.2){$*$};
\draw[shift={(0,0.4)}] (4.8,1.2) node{\tiny $r$};
\node[draw, shape=circle, inner sep=1.8pt] (E) at (6,0){$*$};
\draw[shift={(0,0.4)}] (6,0) node{\tiny $r$};
\node[draw, shape=circle, inner sep=1.8pt] (F) at (4.8,-1.2){$*$};
\draw[shift={(0,0.4)}] (4.8,-1.2) node{\tiny $r$};
\draw node (G) at (3,-1.2){$\cdots$};
\node[draw, shape=circle, inner sep=1.8pt] (H) at (1.2,-1.2){$*$};
\draw[shift={(0,0.4)}] (1.2,-1.2) node{\tiny $r$};

\draw (A) --node[auto=left]{\tiny $\langle r \rangle$} (B);
\draw (B) --node[auto=left]{\tiny $\langle r \rangle$} (C);
\draw (C) --node[auto=left]{\tiny $\langle r \rangle$} (D);
\draw (D) --node[auto=left]{\tiny $\langle r \rangle$} (E);
\draw (A) --node[auto=left]{\tiny $\langle r \rangle$} (H);
\draw (H) --node[auto=left]{\tiny $\langle r \rangle$} (G);
\draw (G) --node[auto=left]{\tiny $\langle r \rangle$} (F);
\draw (F) --node[auto=left]{\tiny $\langle r \rangle$} (E);

\node[draw, shape=circle, inner sep=1.8pt] (B1) at (8,0){$*$};
    \draw[shift={(0,0.4)}] (8,0) node{\tiny $r$};
    \draw[shift={(0,-0.4)}] (8,0) node{\tiny $0$};

    \node[draw, shape=circle, inner sep=1.8pt] (B2) at (10,0){$*$};
    \draw[shift={(0,0.4)}] (10,0) node{\tiny $r$};
    \draw[shift={(0,-0.4)}] (10,0) node{\tiny $0$};

    \draw[shift={(0,0.25)}] (9,0) node{\tiny $\langle 2r \rangle$};
    \draw (B1)--(B2);
    \end{tikzpicture}
    \]
    
On the other hand, we use the term \lq\lq twisted cusp\rq\rq for the following graphs:
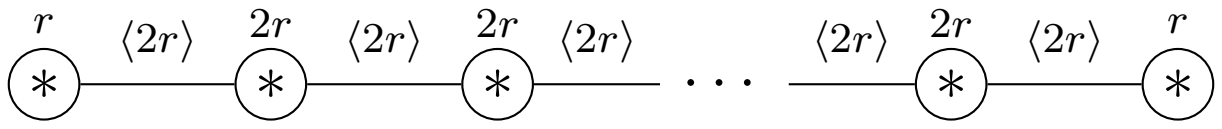
\[
\begin{tikzpicture}
\node[draw, shape=circle, inner sep=1.8pt] (X) at (-1.5,0){$*$};
\draw[shift={(0,0.4)}] (-1.5,0) node{\tiny $r$};
\node[draw, shape=circle, inner sep=1.8pt] (A) at (0,0){$*$};
\draw[shift={(0,0.4)}] (0,0) node{\tiny $2r$};
\node[draw, shape=circle, inner sep=1.8pt] (B) at (1.5,0){$*$};
\draw[shift={(0,0.4)}] (1.5,0) node{\tiny $2r$};
\draw node (C) at (3,0){$\cdots$};
\node[draw, shape=circle, inner sep=1.8pt] (D) at (4.5,0){$*$};
\draw[shift={(0,0.4)}] (4.5,0) node{\tiny $2r$};
\node[draw, shape=circle, inner sep=1.8pt] (E) at (6,0){$*$};
\draw[shift={(0,0.4)}] (6,0) node{\tiny $r$};

\draw (X) --node[auto=left]{\tiny $\langle 2r \rangle$} (A);
\draw (A) --node[auto=left]{\tiny $\langle 2r \rangle$} (B);
\draw (B) --node[auto=left]{\tiny $\langle 2r \rangle$} (C);
\draw (C) --node[auto=left]{\tiny $\langle 2r \rangle$} (D);
\draw (D) --node[auto=left]{\tiny $\langle 2r \rangle$} (E);
\end{tikzpicture}
\]
\end{rem}

Let $\Gamma$ be a graph and let $V$ be the set of vertices of $\Gamma$.
For any pair of distinct vertices $v, w \in V$ ($v \neq w$), we denote by $e_{v,w}$ the number of edges between $v$ and $w$, and by $(r_v, b_v, g_v)$ the weight at $v$.
A bijection $\sigma \colon V \to V$ is called an automorphism of $\Gamma$ if $e_{\sigma(v), \sigma(w)} = e_{v,w}$ for every pair of distinct vertices $v, w \in V$ and $(r_{\sigma(v)}, b_{\sigma(v)}, g_{\sigma(v)}) = (r_v, b_v, g_v)$ for every $v \in V$.
In this case, the \emph{twisted quotient} $\Gamma_{\sigma}$ of $\Gamma$ by $\sigma$ is defined as follows:
\begin{itemize}
    \item The set of vertices of $\Gamma_{\sigma}$ is the quotient set $V / \langle \sigma \rangle$ of $V$ by the cyclic group $\langle \sigma \rangle$.
    \item For an equivalence class $[v] \in V / \langle \sigma \rangle$, the weight at $[v]$ is $(\#[v] \cdot r_v, \#[v] \cdot b_v, g_v)$, where $\#[v]$ denotes the cardinality of the class $[v]$.
    \item For two distinct equivalence classes $\alpha \neq \beta \in V / \langle \sigma \rangle$, the number of edges between them is 
    \[
    e_{\alpha, \beta} \coloneqq \sum_{v,w} e_{v,w}
    \]
    where the sum is taken over all $v,w \in V$ with $[v]=\alpha$ and $[w]=\beta$.
\end{itemize}
According to the classification of dual graphs (\cite[Subsection 3.3]{Kol13} or \cite[Appendix A]{Sato23}), for every lc surface singularity $X$, there exists a unique connected non-twisted graph $\Gamma$ such that the dual graph of $X$ is isomorphic to the twisted quotient of $\Gamma$ by some automorphism $\sigma$.

The following examples illustrate that a twisted graph need not be the graph obtained from the geometric exceptional divisor by Galois descent.

\begin{lem}\label{lem:geom.irr}
We use the notation of Setting \ref{notation:normal surface sings}.
Let $X' = \Spec R' \to X$ be a finite \'etale cover from the spectrum of a local ring $(R', \m', \kappa')$.
We denote by $\Gamma$ (resp.~$\Gamma'$) the dual graph of the minimal resolution of $X$ (resp.~$X'$).
Then the types of $\Gamma$ and $\Gamma'$ are the same, except in the case where $X$ is simple elliptic.
More precisely, we have the following:
\begin{itemize}
    \item If $\Gamma$ is a chain or a twisted chain, then so is $\Gamma'$.
    \item If $\Gamma$ is star-shaped or twisted star-shaped of type $(a,b,c)$ for some $a,b,c \ge 2$, then $\Gamma'$ is of the same type.
    \item If $\Gamma$ is $\tilde{D}_n$ or twisted $\tilde{D}_n$ for some $n \ge 4$, then $\Gamma'$ is $\tilde{D}_m$ or twisted $\tilde{D}_m$ for some $m \ge 4$. Moreover, if $n\geq5$, then $m\geq5$.
    \item If $\Gamma$ is a cusp or a twisted cusp, then so is $\Gamma'$.
\end{itemize}
\end{lem}

\begin{proof}
    The base change $Y\times_X X'\to X'$ is again a minimal resolution since regularity and relative nefness of the canonical divisor are preserved by \'etale base change.
    We first describe the splitting of a twisted chain in detail, and then treat the remaining graph types.
    A chain consisting of a single genus-zero component remains a chain, because a reduced plane conic has at most two components after a separable extension. We now treat a twisted chain with at least two components.
    The dual graph of $\pi$ is as follows:
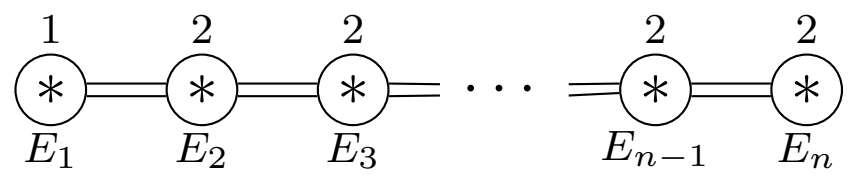
\[
\begin{tikzpicture}
\node[draw, shape=circle, inner sep=1.8pt] (E1) at (1,0){$*$};
\node[draw, shape=circle, inner sep=1.8pt] (E2) at (2,0){$*$};
\node[draw, shape=circle, inner sep=1.8pt] (E3) at (3,0){$*$};
\draw node (dot) at (4,0){$\cdots$};
\node[draw, shape=circle, inner sep=1.8pt] (E_{n-1}) at (5,0){$*$};
\node[draw, shape=circle, inner sep=1.8pt] (E_n) at (6,0){$*$};

\draw[shift={(0,0.4)}] (1,0) node{\tiny $1$};
\draw[shift={(0,0.4)}] (2,0) node{\tiny $2$};
\draw[shift={(0,0.4)}] (3,0) node{\tiny $2$};
\draw[shift={(0,0.4)}] (4,0) node{\tiny ${}$};
\draw[shift={(0,0.4)}] (5,0) node{\tiny $2$};
\draw[shift={(0,0.4)}] (6,0) node{\tiny $2$};

\draw[shift={(0,-0.4)}] (1,0) node{\tiny $E_1$};
\draw[shift={(0,-0.4)}] (2,0) node{\tiny $E_2$};
\draw[shift={(0,-0.4)}] (3,0) node{\tiny $E_3$};
\draw[shift={(0,-0.4)}] (4,0) node{\tiny ${}$};
\draw[shift={(0,-0.4)}] (5,0) node{\tiny $E_{n-1}$};
\draw[shift={(0,-0.4)}] (6,0) node{\tiny $E_n$};

\draw (E1.-10)--(E2.-170);
\draw (E1.10)--(E2.170);
\draw (E2.-10)--(E3.-170);
\draw (E2.10)--(E3.170);
\draw (E3.-10)--(dot.-175);
\draw (E3.10)--(dot.175);
\draw (dot.-5)--(E_{n-1}.-175);
\draw (dot.5)--(E_{n-1}.170);
\draw (E_{n-1}.-10)--(E_n.-170);
\draw (E_{n-1}.10)--(E_n.170);

\end{tikzpicture}
        \]

We set $K \coloneqq H^0(E_2,\sO_{E_2})$.
Since $\sO_{E_1 \cap E_2}$ is a $K$-algebra with dimension two over $\kappa$, we have $E_1 \cap E_2 \cong \Spec K$.
Similarly, we have 
\[
K=H^0(E_2,\sO_{E_2}) \cong \sO_{E_2 \cap E_3} \cong H^0(E_3,\sO_{E_3}) \cong \cdots \cong H^0(E_n,\sO_{E_n}).
\]
In particular, for every $i \ge 2$, since $E_i$ is a genus zero curve over $K$ admitting a Cartier divisor of degree one, we have $E_i \cong \PP^1_K$ \cite[Lemma 10.6]{Kol13}.
On the other hand, it also follows from \cite[Lemma 10.6]{Kol13} that $E_1$ is a degree two curve in $\PP^2_\kappa$.

\textbf{Case 1:} We first consider the case where the base change $(E_1)_{\kappa'} = E_1 \times_{\kappa} \Spec \kappa'$ is irreducible.
If $K \otimes_{\kappa} \kappa'$ is a field, then the base change $(E_i)_{\kappa'} \cong \PP^1_{K \otimes_{\kappa} \kappa'}$ is irreducible for every $i \ge 2$, which implies that we have $\Gamma' = \Gamma$.

From now on, we assume that $K \otimes_{\kappa} \kappa'$ is not a field.
Since $K \otimes_{\kappa} \kappa'$ is a reduced ring whose dimension over $\kappa'$ is two, we have 
\[
K \otimes_{\kappa} \kappa' \cong \kappa' \times \kappa'.
\]
Therefore, one has  
\begin{align*}
    (E_i)_{\kappa'} & \cong E_i^{+} \sqcup E_i^{-}
\end{align*}
for all $i \ge 2$, where $E_i^+ \cong E_i^- \cong \PP^1_{\kappa'}$.
Similarly, we have
\[
(E_i)_{\kappa'} \cap (E_{i+1})_{\kappa'} \cong (E_i \cap E_{i+1})_{\kappa'} \cong\Spec \kappa' \sqcup \Spec \kappa'
\]
for every $i \ge 1$.
After exchanging $E_i^+$ and $E_i^-$ if necessary, we may assume that $E_i^+ \cap E_{i+1}^+ \neq \emptyset$ for every $i \ge 2$.
Then for every $2 \le i < j$ and $a,b \in \{+,-\}$, we have
\begin{align}\label{eq:intersection1}
E_i^a \cap E_j^b \cong \left\{\begin{array}{ll} 
 \Spec \kappa' & \textup{(if $j=i+1$ and $a=b$)} \\
   \emptyset & \textup{(otherwise)}
\end{array}\right. .
\end{align}
Similarly, we have $E_1 \cap E_2^+ \cong E_1 \cap E_2^- \cong \Spec \kappa'$.
The assertion now follows since the dual graph of $X'$ is a chain as follows:
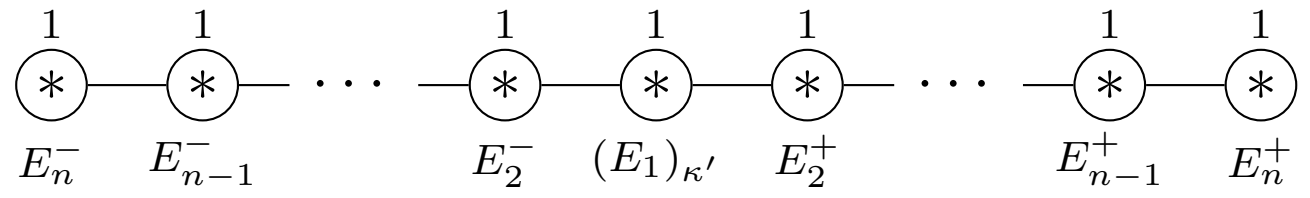
\[
\begin{tikzpicture}
\node[draw, shape=circle, inner sep=1.8pt] (-E_n) at (-2,0){$*$};
\node[draw, shape=circle, inner sep=1.8pt] (-E_{n-1}) at (-1,0){$*$};
\draw node (dot2) at (0,0){$\cdots$};
\node[draw, shape=circle, inner sep=1.8pt] (E1) at (1,0){$*$};
\node[draw, shape=circle, inner sep=1.8pt] (E2) at (2,0){$*$};
\node[draw, shape=circle, inner sep=1.8pt] (E3) at (3,0){$*$};
\draw node (dot) at (4,0){$\cdots$};
\node[draw, shape=circle, inner sep=1.8pt] (E_{n-1}) at (5,0){$*$};
\node[draw, shape=circle, inner sep=1.8pt] (E_n) at (6,0){$*$};

\draw[shift={(0,0.4)}] (-2,0) node{\tiny $1$};
\draw[shift={(0,0.4)}] (-1,0) node{\tiny $1$};
\draw[shift={(0,0.4)}] (1,0) node{\tiny $1$};
\draw[shift={(0,0.4)}] (2,0) node{\tiny $1$};
\draw[shift={(0,0.4)}] (3,0) node{\tiny $1$};
\draw[shift={(0,0.4)}] (4,0) node{\tiny ${}$};
\draw[shift={(0,0.4)}] (5,0) node{\tiny $1$};
\draw[shift={(0,0.4)}] (6,0) node{\tiny $1$};

\draw[shift={(0,-0.5)}] (-2,0) node{\tiny $E_n^-$};
\draw[shift={(0,-0.5)}] (-1,0) node{\tiny $E_{n-1}^-$};
\draw[shift={(0,-0.5)}] (1,0) node{\tiny $E_{2}^-$};
\draw[shift={(0,-0.5)}] (2,0) node{\tiny $(E_1)_{\kappa'}$};
\draw[shift={(0,-0.5)}] (3,0) node{\tiny $E_2^+$};
\draw[shift={(0,-0.5)}] (4,0) node{\tiny ${}$};
\draw[shift={(0,-0.5)}] (5,0) node{\tiny $E_{n-1}^+$};
\draw[shift={(0,-0.5)}] (6,0) node{\tiny $E_n^+$};

\draw (-E_n.0)--(-E_{n-1}.180);
\draw (-E_{n-1}.0)--(dot2.180);
\draw (dot2.0)--(E1.180);
\draw (E1.0)--(E2.-180);
\draw (E2.0)--(E3.-180);
\draw (E3.0)--(dot.-180);
\draw (dot.0)--(E_{n-1}.-180);
\draw (E_{n-1}.0)--(E_n.-180);

\end{tikzpicture}
\]

\textbf{Case 2:} We next consider the case where $(E_1)_{\kappa'}$ is not irreducible.
By Lemma \ref{lem:genus zero curve} below, $E_1$ has a single singular point $Q \in E_1$, whose residue field is $\kappa$, $(E_1)_{\kappa'}$ is isomorphic to a union of two lines $E_1^+$ and $E_1^-$ in $\PP^2_{\kappa'}$, and $E_1^+ \cap E_1^-$ consists of the single point $Q'$ lying over $Q$.

Since we have $E_1 \cap E_2 \cong \Spec K \not\cong \Spec \kappa(Q)$, it also follows from Lemma \ref{lem:genus zero curve} that $(E_1)_{\kappa'} \cap (E_2)_{\kappa'}$ consists of two points $P^+$ and $P^-$, which are not equal to $Q'$.
Noting that $(E_1)_{\kappa'} \cap (E_2)_{\kappa'} \cong \Spec K \otimes_{\kappa} \kappa'$, the base change $K \otimes_{\kappa} \kappa'$ is not a field.
Therefore, as in Case 1, for each $i \ge 2$, the base change $(E_i)_{\kappa'}$ is a disjoint union of two lines $E_i^+$ and $E_i^-$.

By swapping $P^+$ with $P^-$, and $E_i^+$ with $E_i^-$ if necessary, we conclude that the dual graph of $X'$ is a chain given as follows:
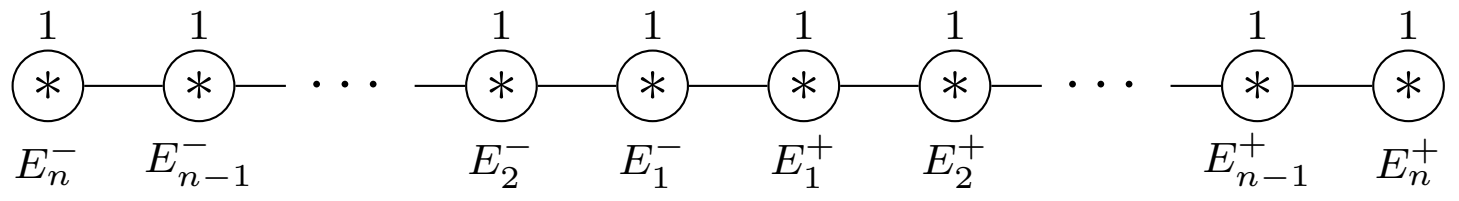
\[
\begin{tikzpicture}
\node[draw, shape=circle, inner sep=1.8pt] (-E_n) at (-3,0){$*$};
\node[draw, shape=circle, inner sep=1.8pt] (-E_{n-1}) at (-2,0){$*$};
\draw node (dot2) at (-1,0){$\cdots$};
\node[draw, shape=circle, inner sep=1.8pt] (E0) at (0,0){$*$};
\node[draw, shape=circle, inner sep=1.8pt] (E1) at (1,0){$*$};
\node[draw, shape=circle, inner sep=1.8pt] (E2) at (2,0){$*$};
\node[draw, shape=circle, inner sep=1.8pt] (E3) at (3,0){$*$};
\draw node (dot) at (4,0){$\cdots$};
\node[draw, shape=circle, inner sep=1.8pt] (E_{n-1}) at (5,0){$*$};
\node[draw, shape=circle, inner sep=1.8pt] (E_n) at (6,0){$*$};

\draw[shift={(0,0.4)}] (-3,0) node{\tiny $1$};
\draw[shift={(0,0.4)}] (-2,0) node{\tiny $1$};
\draw[shift={(0,0.4)}] (0,0) node{\tiny $1$};
\draw[shift={(0,0.4)}] (1,0) node{\tiny $1$};
\draw[shift={(0,0.4)}] (2,0) node{\tiny $1$};
\draw[shift={(0,0.4)}] (3,0) node{\tiny $1$};
\draw[shift={(0,0.4)}] (4,0) node{\tiny ${}$};
\draw[shift={(0,0.4)}] (5,0) node{\tiny $1$};
\draw[shift={(0,0.4)}] (6,0) node{\tiny $1$};

\draw[shift={(0,-0.5)}] (-3,0) node{\tiny $E_n^-$};
\draw[shift={(0,-0.5)}] (-2,0) node{\tiny $E_{n-1}^-$};
\draw[shift={(0,-0.5)}] (0,0) node{\tiny $E_{2}^-$};
\draw[shift={(0,-0.5)}] (1,0) node{\tiny $E_1^-$};
\draw[shift={(0,-0.5)}] (2,0) node{\tiny $E_1^+$};
\draw[shift={(0,-0.5)}] (3,0) node{\tiny $E_2^+$};
\draw[shift={(0,-0.5)}] (4,0) node{\tiny ${}$};
\draw[shift={(0,-0.5)}] (5,0) node{\tiny $E_{n-1}^+$};
\draw[shift={(0,-0.5)}] (6,0) node{\tiny $E_n^+$};

\draw (-E_n.0)--(-E_{n-1}.180);
\draw (-E_{n-1}.0)--(dot2.180);
\draw (dot2.0)--(E0.180);
\draw (E0.0)--(E1.180);
\draw (E1.0)--(E2.-180);
\draw (E2.0)--(E3.-180);
\draw (E3.0)--(dot.-180);
\draw (dot.0)--(E_{n-1}.-180);
\draw (E_{n-1}.0)--(E_n.-180);

\end{tikzpicture}
\]
It remains to discuss the other graph types.
For a non-twisted chain with at least two vertices, each component has a Cartier divisor of degree one over $\kappa$ and is therefore $\PP^1_\kappa$; its graph is unchanged by the extension.

For a star-shaped graph or one of its twisted forms, the central component has a Cartier divisor of odd-degree: the classified graphs provide an attachment of degree one, except for the forms with a single triple attachment, where the degree is three.
Thus the central component is $\PP^1_\kappa$ by \cite[Lemma 10.6]{Kol13}.
Every component in an arm is $\PP^1$ over its field of constants, since an adjacent intersection gives a Cartier divisor of degree one over that field.
Under the separable extension $\kappa'/\kappa$, a field of constants becomes a product of fields.
The corresponding arm therefore splits into disjoint copies of the same arm, with the same self-intersection numbers normalized by the degrees of the constant fields.
The central component remains irreducible, and the attachments split accordingly.
Hence the three arm types, including their determinants, are unchanged.
This applies equally to the spherical types in \cite[Figure 1]{Sato23} and the parabolic types in \cite[Figures 2 and 3]{Sato23}.

For the $\tilde D$ types, write $K_Y+\Delta_Y=\pi^*K_X$.
Their exceptional coefficients in $\Delta_Y$ are all $1/2$ or $1$, with at least one coefficient $1/2$.
These coefficients are unchanged by \'etale base change.
Moreover, rationality is preserved, since $R^1\pi_*\sO_Y$ commutes with the flat base change $X'\to X$.
Among the rational non-klt graphs in \cite[Figures 2 and 3]{Sato23}, the $\tilde D$ types are precisely those with this coefficient pattern: each parabolic star has an arm coefficient different from $1/2$ and $1$.
Thus $\Gamma'$ is again of $\tilde D$ type.

Finally, a cusp or twisted cusp is nonrational, and each exceptional component has vanishing $H^1$.
Both properties persist after the present base change.
Indeed, $(E_i)_{\kappa'}$ is reduced, and the exact sequence comparing its structure sheaf with those of its irreducible components shows that $H^1((E_i)_{\kappa'},\sO)=0$ forces the $H^1$ of every component to vanish.
The resulting graph is therefore a nonrational lc graph all of whose vertices have genus zero.
By \cite[Figure 4]{Sato23}, it is a cusp or a twisted cusp, rather than a simple elliptic graph.

We explain also why the cases $n\geq5$ cannot become $\tilde D_4$ or a twisted $\tilde D_4$.
Write $K_Y+\Delta_Y=\pi^*K_X$ and set
\[
w(X):=\sum_{\operatorname{coeff}_{E_i}\Delta_Y=1}
       [H^0(E_i,\sO_{E_i}):\kappa].
\]
Inspection of the $\tilde D_n$ graphs and their twisted forms in \cite[Figures 2 and 3]{Sato23} shows that $w(X)=1$ for the $\tilde D_4$ forms, whereas $w(X)\geq2$ for the forms with $n\geq5$.
Under a finite separable extension $\kappa'/\kappa$, the numerical discrepancies pull back unchanged.
For each $E_i$, the reduced base change $(E_i)_{\kappa'}$ has irreducible components $E_{ij}'$, and restriction of functions gives an injection
\[
H^0(E_i,\sO_{E_i})\otimes_\kappa\kappa'
\hookrightarrow\prod_j H^0(E_{ij}',\sO_{E_{ij}'}).
\]
Thus $w(X')\geq w(X)$, proving the assertion.
\end{proof}

\begin{lem}\label{lem:genus zero curve}
Let $C$ be an integral Gorenstein projective curve over a field $K=H^0(C,\sO_C)$ of genus zero.
We further assume that $C$ is geometrically reduced over $K$, but not isomorphic to $\PP^1_K$.
\begin{enumerate}[label=\textup{(\arabic*)}]
    \item There exists a smooth closed point $P \in C$ such that the residue field $L \coloneqq \kappa(P)$ is a separable extension over $K$ of degree two.
    \item If $C$ is geometrically integral over $K$, the base change $C_L \coloneqq C \times_K \Spec L$ is isomorphic to $\PP^1_L$.
    \item If $C$ is not geometrically irreducible over $K$, then $C_L$ is not irreducible.
    Moreover, we have the following properties:
    \begin{enumerate}[label=\textup{(\roman*)}]
    \item $C_L = C_L^+ \cup C_L^-$ with $C_L^+ \cong C_L^- \cong \PP^1_L$ and $C_L^+ \cap C_L^- \cong \Spec L$.
    \item There exists an isomorphism $\sigma: C_L^+ \xrightarrow{\sim} C_L^-$ over $C$.
    \item For any point $Q \in C_L^+$, we have $\nu^{-1}(\nu(Q))=\{Q, \sigma(Q)\}$, where $\nu: C_L \to C$ is the natural morphism.
    \end{enumerate}
\end{enumerate}
\end{lem}

\begin{proof}
    It follows from \cite[Lemma 10.6]{Kol13} that $C$ is isomorphic to a degree two curve in $\PP^2_{K}$.
    We first prove (1).
    
    \textbf{Case 1:}
    We first consider the case where the characteristic of $K$ is not $2$.
    After a change of variables, we can write
    \[
    C \cong Z(x^2+ay^2+bz^2)
    \]
    for some $a,b \in K$.
    Noting that $C$ is non-reduced if $a=b=0$, we may assume that $a \neq 0$.
    By the Jacobian criterion, $C$ is smooth if $b \neq 0$ and smooth outside $[0,0,1]$ if $b=0$.
    Since a smooth $K$-rational point would give an invertible sheaf of degree one and hence an isomorphism $C\cong\PP^1_K$ \cite[Lemma 10.6]{Kol13}, there is no such point.
    This shows that $x^2+ay^2$ is an irreducible polynomial in $K[x,y]$ and the closed point $P \coloneqq Z(x^2+ay^2, z)$ satisfies the assertion in (1).

    \textbf{Case 2:}
    We next consider the case where the characteristic is $2$.
    We write 
    \[
    C \cong Z(ax^2+by^2+cz^2+dxy+eyz+fzx)
    \]
    for some $a,b,c,d,e,f \in K$.
    Since $C$ is geometrically reduced over $K$, at least one of $d,e$ and $f$ is non-zero.
    We may assume that $d=1$.
    After a change of variables, we can write
    \[
    C \cong Z(ax^2+by^2+cz^2+xy)
    \]
    for some $a,b,c \in K$.
    By the Jacobian criterion, $C$ is smooth if $c \neq 0$ and smooth outside $[0,0,1]$ if $c=0$.
    Since a smooth $K$-rational point would give an invertible sheaf of degree one and hence an isomorphism $C\cong\PP^1_K$ \cite[Lemma 10.6]{Kol13}, there is no such point.
    This shows that $ax^2+by^2+xy$ is an irreducible polynomial in $K[x,y]$ and the closed point $P \coloneqq Z(ax^2+by^2+xy, z)$ satisfies the assertion in (1).

    Since we have $L \otimes_K L \cong L \times L$, the pullback $\nu^{-1}(P)$ of $P$ by the natural morphism $\nu : C_L \to C$ is a set of two smooth $L$-rational points.
    In particular, if $C_L$ is irreducible, then $C_L$ is isomorphic to $\PP^1_L$ by \cite[Lemma 10.6]{Kol13}.
    This shows the assertion in (2).
    
    For (3), we assume that $C$ is not geometrically irreducible over $K$. 
    If $C_L$ is irreducible, then by the above argument, $C_L$ is isomorphic to $\PP^1_L$, which contradicts the fact that $C_L$ is not geometrically irreducible over $L$.
    Therefore, $C_L$ is not irreducible.
    The assertion in (3-i) follows from the fact that $C_L$ is a reduced reducible Cartier divisor of degree two in $\PP^2_L$.

    For (3-ii), we consider the generator $g \colon \Spec L \xrightarrow{\sim} \Spec L$ of the Galois group $\Z/2\Z$ of $L$ over $K$ and we write
    \[
    \sigma \coloneqq \id_C \times_K g : C_L \xrightarrow{\sim} C_L.
    \]
    It suffices to show that $\sigma(\eta^+)=\eta^-$, where $\eta^+$ (resp.~$\eta^-$) is the generic point of $C_L^+$ (resp.~$C_L^-$).
    Let $K(C)$ be the function field of $C$.
    Since we have $K(C) \otimes_K L \cong K(C) \times K(C)$, the isomorphism $\sigma$ induces the $K(C)$-algebra automorphism 
    \[
    \sigma^* : K(C) \times K(C) \xrightarrow{\sim} K(C) \times K(C).
    \]
    We note that a $K(C)$-algebra automorphism on $K(C) \times K(C)$ is either the identity map or the map defined by $(a,b) \mapsto (b,a)$.
    Noting that $\sigma^*$ is not the identity map, $\sigma^*$ coincides with the latter.
    This implies that $\sigma(\eta^+)=\eta^-$.
    
    For (3-iii), we first note that for any field $\kappa$ containing $K$, the base change $\kappa \otimes_K L$ is either a field extension over $\kappa$ of degree two or isomorphic to $\kappa \times \kappa$.
    This shows that the fiber $\nu^{-1}(R)$ of any point $R \in C$ contains at most two points.
    Take a point $Q \in C_L^+$.
    We note that the set $\{Q, \sigma(Q)\}$ is contained in $\nu^{-1}(\nu(Q))$.
    Assume to the contrary that there is a point $Q' \in \nu^{-1}(\nu(Q)) \setminus \{Q, \sigma(Q)\}$.
    Then we have the inclusion
    \[
    \{Q, \sigma(Q), Q', \sigma(Q')\} \subseteq \nu^{-1}(\nu(Q)).
    \]
    Combining this with the fact that $\# \nu^{-1}(\nu(Q)) \le 2$, we have $Q=\sigma(Q)$ and $Q' = \sigma(Q')$.
    Since $Q \in C_L^+$ and $\sigma(Q) \in \sigma(C_L^+)=C_L^-$, we have $Q \in C_L^+ \cap C_L^-$.
    Similarly, one has $Q' \in C_L^+ \cap C_L^-$.
    This contradicts the fact that $C_L^{+} \cap C_L^-$ consists of a single point.
    \end{proof}

\subsubsection{Rational double points}
We use the notation of Setting \ref{notation:normal surface sings}.
By \cite[Proposition 3.7]{Sato-Takagi(generalhyperplane)}, we have
\[
\dim_{\kappa} H^0(E_i, \sO_{E_i})-\dim_{\kappa} H^1(E_i, \sO_{E_i})=\chi(E_i, \sO_{E_i})=-\frac{E_i^2+K_Y\cdot E_i}{2}.
\]
Since $E_i^2<0$, it follows that $K_Y\cdot E_i=0$ if and only if $E_i^2=-2\dim_{\kappa} H^0(E_i, \sO_{E_i})$ and $H^1(E_i, \sO_{E_i})=0$.
We say that $E_i$ is a \textit{$(-2)$-curve} if these equivalent conditions hold.
We also say that $X$ is a \emph{rational double point} (\emph{RDP}, for short) if each $E_i$ is a $(-2)$-curve.

\begin{prop}\label{prop:RDP}
    Suppose that $X$ is klt.
    Then the following are equivalent:
    \begin{enumerate}
        \item[\textup{(1)}] $X$ is an RDP.
        \item[\textup{(2)}] $X$ is canonical.
        \item[\textup{(3)}] $X$ is Gorenstein.
    \end{enumerate}
\end{prop}
\begin{proof}
Because $X$ is klt, if $K_X$ is Cartier, then all discrepancies are nonnegative integers. Thus (3) implies (2).
The implication (2)$\Rightarrow$(3) is \cite[Proposition 3.9]{Sato-Takagi(generalhyperplane)}.

For the remaining equivalence, write $K_Y=\pi^*K_X+A$ with $A$ exceptional.
It then follows from the negative definiteness of the exceptional intersection that $X$ is canonical if and only if $A = 0$.
This is equivalent to $K_Y\cdot E_i=0$ for every $i$, that is, each $E_i$ is a $(-2)$-curve.
\end{proof}

\section{Geometric reducedness of exceptional curves}\label{sec:geometric-reducedness}
We study geometric reducedness over the residue field of exceptional curves on minimal resolutions. We first observe the possible failures for lc surface singularities, and then show how $F$-purity and strong $F$-regularity exclude these failures.

\subsection{Geometric reducedness and separability}\label{subsec:geom red}

We begin by relating geometric reducedness to the separability of the field of global functions.

\begin{lem}\label{lem:geom.red to separable}
    Let $X$ be a proper variety over a field $k$.
    If $X$ is geometrically reduced over $k$, then $H^0(X, \sO_X)$ is a finite separable extension of $k$.
\end{lem}
\begin{proof}
    Set $L\coloneqq H^0(X,\sO_X)$, which is a finite field extension of $k$ because $X$ is proper and integral.
    For every field extension $K/k$, flat base change gives
    \[
    L\otimes_k K\cong H^0(X\times_k K,\sO_{X\times_k K}).
    \]
    The right-hand side is reduced by assumption. Thus $L/k$ is separable (see \cite[Section 26]{Matsumura}).
\end{proof}

For the exceptional curves of an lc surface singularity, a converse holds outside a short list of cases:

\begin{prop}\label{prop:equivalence of geom.red and separability}
    We use the notation of Setting \ref{notation:normal surface sings}.
    If each field $K_i \coloneqq H^0(E_i,\mathcal{O}_{E_i})$ is separable over $\kappa$, then every $E_i$ is geometrically reduced over $\kappa$, except when one of the following holds:
    \begin{enumerate}[label=\textup{(\alph*)}]
        \item $X$ is simple elliptic with $p \le 3$,
        \item $X$ is a cusp of length $2$ with $p=2$, or
        \item the dual graph of $\pi$ is a chain of length $1$ with $p=2$.
    \end{enumerate}
    \begin{figure}[h]
\centering
\begin{minipage}[b]{0.32\columnwidth}
    \centering
        \[
    \begin{tikzpicture}
        \node[draw, shape=circle, inner sep=1.8pt] (A1) at (0,0){$*$};
    \draw[shift={(0,0.4)}] (0,0) node{\tiny $1$};
    \draw[shift={(0,-0.4)}] (0,0) node{\tiny $1$};
    \end{tikzpicture}
    \]
\end{minipage}
\begin{minipage}[b]{0.32\columnwidth}
    \centering
        \[
    \begin{tikzpicture}
        \node[draw, shape=circle, inner sep=1.8pt] (A1) at (0,0){$*$};
    \draw[shift={(0,0.4)}] (0,0) node{\tiny $1$};
    \draw[shift={(0,-0.4)}] (0,0) node{\tiny $0$};
    \end{tikzpicture}
    \]
\end{minipage}
\begin{minipage}[b]{0.32\columnwidth}
    \centering
    \[
    \begin{tikzpicture}
    \node[draw, shape=circle, inner sep=1.8pt] (B1) at (0,0){$*$};
    \draw[shift={(0,0.4)}] (0,0) node{\tiny $r$};
    \draw[shift={(0,-0.4)}] (0,0) node{\tiny $0$};

    \node[draw, shape=circle, inner sep=1.8pt] (B2) at (2,0){$*$};
    \draw[shift={(0,0.4)}] (2,0) node{\tiny $r$};
    \draw[shift={(0,-0.4)}] (2,0) node{\tiny $0$};

    \draw[shift={(0,0.25)}] (1,0) node{\tiny $\langle 2r \rangle$};
    \draw (B1)--(B2);
    \end{tikzpicture}
    \]
\end{minipage}
\caption{The three possible exceptions in Proposition \ref{prop:equivalence of geom.red and separability}: simple elliptic, a chain of length one, and a cusp of length two.}
\end{figure}
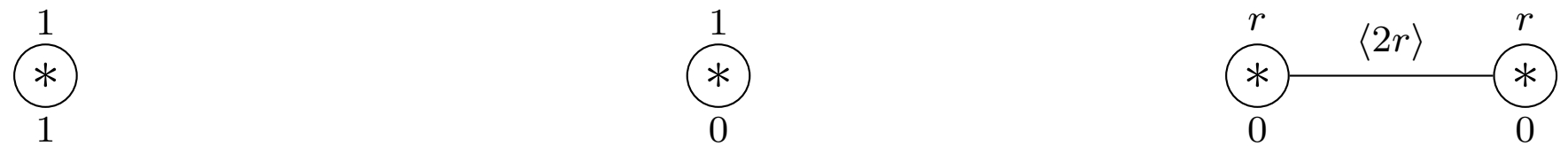
    
\end{prop}

\begin{proof}
    Suppose first that $X$ is simple elliptic and $p>3$, and write $E$ for its exceptional curve.
    If $E$ is regular, the assertion follows from \cite[Theorem 1.1]{PW} (see also \cite[Proposition 9.11(2)]{Tanaka(invariants)}).
    Otherwise, the lc condition implies that the non-regular points of $E$ are nodes, and \cite[Corollary 3.10]{Sato23} gives the assertion.

    Suppose next that $X$ is not simple elliptic.
    By the classification of dual graphs \cite[Appendix A]{Sato23}, every $E_i$ has genus zero.
    Since $Y$ is regular and $E_i$ is a Cartier divisor on $Y$, the curve $E_i$ is Gorenstein.
    If $p\ne 2$, the assertion therefore follows from \cite[Lemma 3.3]{Sato23}.

    It remains to treat $p=2$, excluding cusps of length two and chains of length one.
    The classification shows that for each $E_i$ there is a distinct component $E_j$ with
    \[
    E_i\cdot E_j=\max\{r_i,r_j\}.
    \]
    The coordinate ring of the scheme-theoretic intersection $D\coloneqq E_i\cap E_j$ is an algebra over both $K_i$ and $K_j$.
    Its dimension over $\kappa$ is $\max\{r_i,r_j\}$, so $D$ is isomorphic to the spectrum of whichever of these two fields has larger degree.
    Thus $D$ is geometrically reduced over $\kappa$ by the separability assumption.
    Moreover, $D$ is a nonempty effective Cartier divisor on $E_i$.
    Lemma \ref{lem:geom.red} below now shows that $E_i$ is geometrically reduced over $\kappa$.
\end{proof}

\begin{lem}\label{lem:geom.red}
    Let $Z$ be an integral locally Noetherian scheme over a field $k$.
    Suppose that $Z$ contains a nonempty effective Cartier divisor $D$ that is geometrically reduced over $k$.
    Then $Z$ is geometrically reduced over $k$.
\end{lem}
\begin{proof}
    Let $K/k$ be a finite field extension, choose $P\in\Supp(D)$, and set
    \[
    A\coloneqq\sO_{Z,P},\qquad B\coloneqq A\otimes_k K.
    \]
    Let $f\in A$ be a local equation of $D$.
    The ring $B$ is finite and free over $A$, so it is Noetherian, $f$ is a non-zero-divisor on $B$, and $f$ belongs to the Jacobson radical of $B$.
    By hypothesis, $B/fB\cong(A/fA)\otimes_k K$ is reduced.
    If $b\in B$ is nilpotent, then $b=fb_1$ for some $b_1\in B$ since $f$ is a non-zero-divisor, $b_1$ is also nilpotent.
    Iterating gives $b\in\bigcap_{n\ge1}f^nB=0$ by the Krull intersection theorem.
    Thus $B$ is reduced.

    Writing $k(Z)$ for the function field of $Z$, we conclude that $k(Z)\otimes_k K$ is reduced, since it is a localization of $B$.
    For every nonempty affine open $U\subseteq Z$, the injection
    \[
    \sO_Z(U)\otimes_k K\hookrightarrow k(Z)\otimes_k K
    \]
    shows that $U\times_k K$ is reduced.
    Hence $Z\times_k K$ is reduced for every finite extension $K/k$, as required.
\end{proof}

\begin{prop}\label{prop:geom red,cusp}
    We use the notation of Setting \ref{notation:normal surface sings}.
    \begin{enumerate}[label=\textup{(\arabic*)}]
    \item Assume that $X$ is a cusp of length at least three.
    
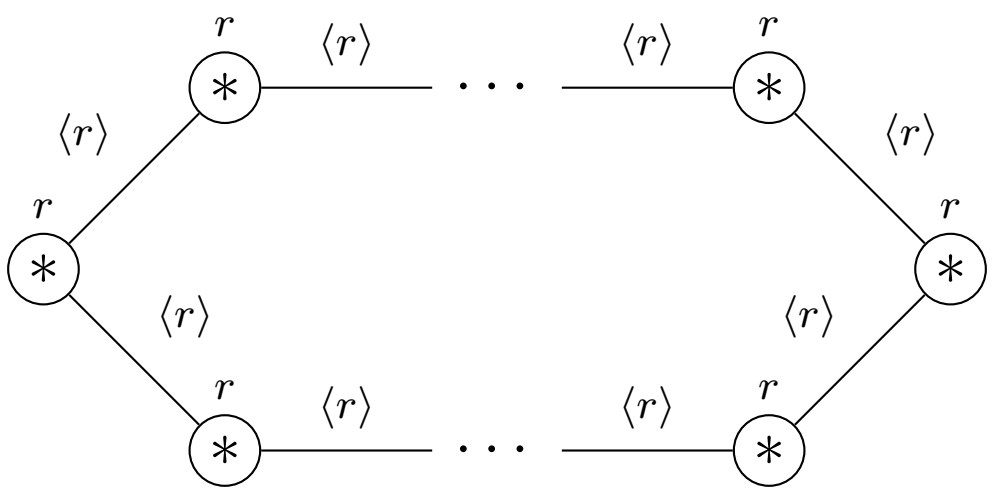
\[
    \begin{tikzpicture}
\node[draw, shape=circle, inner sep=1.8pt] (A) at (0,0){$*$};
\draw[shift={(0,0.4)}] (0,0) node{\tiny $r$};
\node[draw, shape=circle, inner sep=1.8pt] (B) at (1.2,1.2){$*$};
\draw[shift={(0,0.4)}] (1.2,1.2) node{\tiny $r$};
\draw node (C) at (3,1.2){$\cdots$};
\node[draw, shape=circle, inner sep=1.8pt] (D) at (4.8 ,1.2){$*$};
\draw[shift={(0,0.4)}] (4.8,1.2) node{\tiny $r$};
\node[draw, shape=circle, inner sep=1.8pt] (E) at (6,0){$*$};
\draw[shift={(0,0.4)}] (6,0) node{\tiny $r$};
\node[draw, shape=circle, inner sep=1.8pt] (F) at (4.8,-1.2){$*$};
\draw[shift={(0,0.4)}] (4.8,-1.2) node{\tiny $r$};
\draw node (G) at (3,-1.2){$\cdots$};
\node[draw, shape=circle, inner sep=1.8pt] (H) at (1.2,-1.2){$*$};
\draw[shift={(0,0.4)}] (1.2,-1.2) node{\tiny $r$};

\draw (A) --node[auto=left]{\tiny $\langle r \rangle$} (B);
\draw (B) --node[auto=left]{\tiny $\langle r \rangle$} (C);
\draw (C) --node[auto=left]{\tiny $\langle r \rangle$} (D);
\draw (D) --node[auto=left]{\tiny $\langle r \rangle$} (E);
\draw (A) --node[auto=left]{\tiny $\langle r \rangle$} (H);
\draw (H) --node[auto=left]{\tiny $\langle r \rangle$} (G);
\draw (G) --node[auto=left]{\tiny $\langle r \rangle$} (F);
\draw (F) --node[auto=left]{\tiny $\langle r \rangle$} (E);
\end{tikzpicture}
\]
    Then the fields $K_i=H^0(E_i,\sO_{E_i})$ satisfy 
    \[
    K_1 \cong K_2 \cong \dots \cong K_n,
    \]
    and all $K_i$ are separable over $\kappa$.
    
\item Assume that $X$ is a twisted cusp of length $n\ge3$.
    
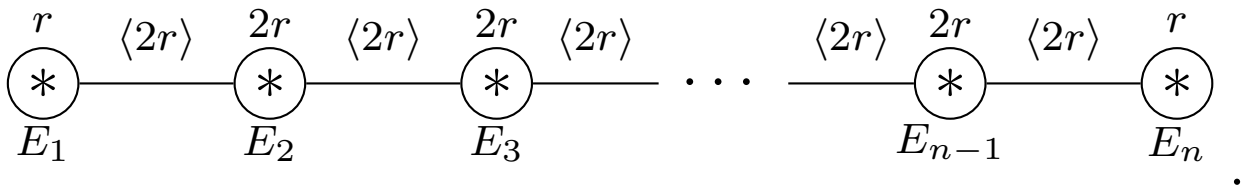
\[
    \begin{tikzpicture}
\node[draw, shape=circle, inner sep=1.8pt] (X) at (-1.5,0){$*$};
\draw[shift={(0,0.4)}] (-1.5,0) node{\tiny $r$};
\draw[shift={(0,-0.4)}] (-1.5,0) node{\tiny $E_1$};
\node[draw, shape=circle, inner sep=1.8pt] (A) at (0,0){$*$};
\draw[shift={(0,0.4)}] (0,0) node{\tiny $2r$};
\draw[shift={(0,-0.4)}] (0,0) node{\tiny $E_2$};
\node[draw, shape=circle, inner sep=1.8pt] (B) at (1.5,0){$*$};
\draw[shift={(0,0.4)}] (1.5,0) node{\tiny $2r$};
\draw[shift={(0,-0.4)}] (1.5,0) node{\tiny $E_3$};
\draw node (C) at (3,0){$\cdots$};
\node[draw, shape=circle, inner sep=1.8pt] (D) at (4.5,0){$*$};
\draw[shift={(0,0.4)}] (4.5,0) node{\tiny $2r$};
\draw[shift={(0,-0.4)}] (4.5,0) node{\tiny $E_{n-1}$};
\node[draw, shape=circle, inner sep=1.8pt] (E) at (6,0){$*$};
\draw[shift={(0,0.4)}] (6,0) node{\tiny $r$};
\draw[shift={(0,-0.4)}] (6,0) node{\tiny $E_n$};

\draw (X) --node[auto=left]{\tiny $\langle 2r \rangle$} (A);
\draw (A) --node[auto=left]{\tiny $\langle 2r \rangle$} (B);
\draw (B) --node[auto=left]{\tiny $\langle 2r \rangle$} (C);
\draw (C) --node[auto=left]{\tiny $\langle 2r \rangle$} (D);
\draw (D) --node[auto=left]{\tiny $\langle 2r \rangle$} (E);
\end{tikzpicture}.
    \]
    Then the intersections of adjacent exceptional curves identify $K_2\cong K_3\cong\cdots\cong K_{n-1}$ and define embeddings of $K_1$ and $K_n$ into this common field.
    If we further assume that $K_2$ is separable over both $K_1$ and $K_n$, then $K_i$ is separable over $\kappa$ for every $i$.
    \item Assume that $X$ is a cusp of length two.
        
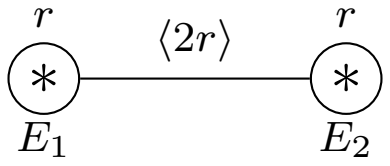
\[
    \begin{tikzpicture}
    \node[draw, shape=circle, inner sep=1.8pt] (B1) at (0,0){$*$};
    \draw[shift={(0,0.4)}] (0,0) node{\tiny $r$};
    \draw[shift={(0,-0.4)}] (0,0) node{\tiny $E_1$};

    \node[draw, shape=circle, inner sep=1.8pt] (B2) at (2,0){$*$};
    \draw[shift={(0,0.4)}] (2,0) node{\tiny $r$};
    \draw[shift={(0,-0.4)}] (2,0) node{\tiny $E_2$};

    \draw[shift={(0,0.25)}] (1,0) node{\tiny $\langle 2r \rangle$};
    \draw (B1)--(B2);
    \end{tikzpicture}
    \]
    Then one of the following holds:
    \begin{itemize}
        \item[\textup{(3-1)}] We have 
        \[
        E_1 \cap E_2 \cong \Spec K_1 \coprod \Spec K_1,
        \] and $K_1 \cong K_2$, which are separable over $\kappa$.
        \item[\textup{(3-2)}] We have $E_1 \cap E_2 \cong \Spec L$ for some field extension $L/K_1$ of degree $2$.
        In this case, if we further assume that $L$ is separable over both $K_1$ and $K_2$, then $K_1$ and $K_2$ are separable over $\kappa$.
    \end{itemize}
    \end{enumerate}
\end{prop}

\begin{proof}
    Assertion (1) is \cite[Proposition 4.8]{Sato23}.
    For (2) and (3), set $E\coloneqq\sum_iE_i$.
    We use $H^0(E,\sO_E)=\kappa$ \cite[Corollary 10.10]{Kol13}.

    In (2), every adjacent intersection has coordinate ring of dimension $2r$ over $\kappa$.
    It is therefore the spectrum of the constant field of the adjacent interior component, whose degree is $2r$.
    This identifies all the interior fields with a common field $L$ and embeds $K_1$ and $K_n$ into $L$ as subfields of index two.
    Gluing global functions along the chain \cite[Lemma 4.7]{Sato23} gives
    \[
    \kappa=H^0(E,\sO_E)=K_1\cap K_n\subseteq L.
    \]
    If $L/K_1$ and $L/K_n$ are separable, they are quadratic Galois extensions.
    Let $\sigma_1$ and $\sigma_n$ be their nontrivial automorphisms, and let $G$ be the subgroup of $\operatorname{Aut}_{\kappa}(L)$ that they generate.
    Since $L/\kappa$ is finite, $G$ is finite, and
    \[
    L^G=L^{\sigma_1}\cap L^{\sigma_n}=K_1\cap K_n=\kappa.
    \]
    The fixed-field theorem shows that $L/\kappa$ is Galois. In particular, every $K_i/\kappa$ is separable.

    For (3), the intersection $E_1\cap E_2$ is reduced: otherwise, the reduced boundary $E_1+E_2$ would not have a node at its intersection point, contradicting the lc condition on $(Y,E_1+E_2)$.
    This is the argument in the proof of \cite[Proposition 4.11]{Sato23}, and it does not use the characteristic assumption in that proposition.
    Since the intersection has degree $2r$, whereas each $K_i$ has degree $r$, it consists either of two points with residue field isomorphic to each $K_i$, or of one point with residue field $L$ satisfying $[L:K_i]=2$.

    In the two-point case, use the restriction maps at the first point to identify $K_1$ with $K_2$.
    Agreement at the second point then identifies $H^0(E,\sO_E)$ with the fixed field of a $\kappa$-automorphism $\sigma$ of $K_1$.
    Thus $K_1^{\langle\sigma\rangle}=\kappa$, so $K_1/\kappa$ is cyclic Galois, proving (3-1).
    In the one-point case, gluing gives $H^0(E,\sO_E)=K_1\cap K_2$ inside $L$.
    If both $L/K_1$ and $L/K_2$ are separable, the same fixed-field argument as in (2) proves that $L/\kappa$ is Galois, and hence that both $K_i/\kappa$ are separable.
\end{proof}

\begin{cor}\label{cor:geom red,p>3}
    In Setting \ref{notation:normal surface sings}, assume $p>3$.
    Then $E_i$ is geometrically reduced over $\kappa$ for all $i$.
\end{cor}

\begin{proof}
    By Proposition \ref{prop:equivalence of geom.red and separability}, it suffices to show that $K_i = H^0(E_i,\sO_{E_i})$ is separable over $\kappa$ for all $i$.
    If $X$ is a cusp or a twisted cusp, this follows from Proposition \ref{prop:geom red,cusp}, since the quadratic extensions appearing there are separable when $p>3$.
    For the remaining cases, it follows from the classification of dual graphs that $r_i=\dim_{\kappa}K_i \le 4$ for every $i$.
    Since $r_i$ is not divisible by $p$, $K_i$ is separable over $\kappa$, as desired.
\end{proof}

\begin{cor}\label{cor:geom red,p=3}    
In Setting \ref{notation:normal surface sings}, assume $p=3$.
Then every $E_i$ is geometrically reduced over $\kappa$, except when one of the following holds:
    \begin{itemize}
        \item $X$ is simple elliptic, or 
        \item the dual graph of $\pi$ is 
    
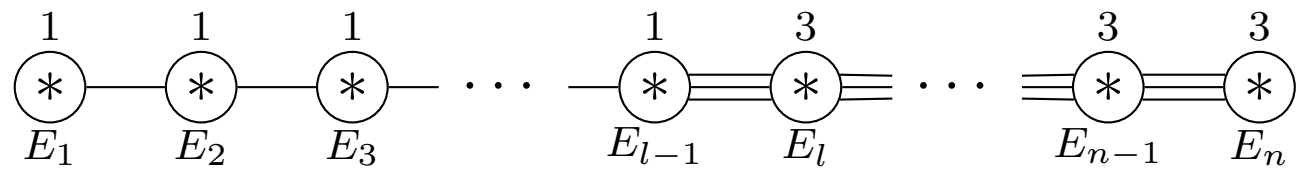
\[
\begin{tikzpicture}
\node[draw, shape=circle, inner sep=1.8pt] (E1) at (1,0){$*$};
\node[draw, shape=circle, inner sep=1.8pt] (E2) at (2,0){$*$};
\node[draw, shape=circle, inner sep=1.8pt] (E3) at (3,0){$*$};
\draw node (dot) at (4,0){$\cdots$};
\node[draw, shape=circle, inner sep=1.8pt] (E_{l-1}) at (5,0){$*$};
\node[draw, shape=circle, inner sep=1.8pt] (E_l) at (6,0){$*$};
\draw node (dot2) at (7,0){$\cdots$};
\node[draw, shape=circle, inner sep=1.8pt] (E_{n-1}) at (8,0){$*$};
\node[draw, shape=circle, inner sep=1.8pt] (E_n) at (9,0){$*$};

\draw[shift={(0,0.4)}] (1,0) node{\tiny $1$};
\draw[shift={(0,0.4)}] (2,0) node{\tiny $1$};
\draw[shift={(0,0.4)}] (3,0) node{\tiny $1$};
\draw[shift={(0,0.4)}] (4,0) node{\tiny ${}$};
\draw[shift={(0,0.4)}] (5,0) node{\tiny $1$};
\draw[shift={(0,0.4)}] (6,0) node{\tiny $3$};
\draw[shift={(0,0.4)}] (8,0) node{\tiny $3$};
\draw[shift={(0,0.4)}] (9,0) node{\tiny $3$};

\draw[shift={(0,-0.4)}] (1,0) node{\tiny $E_1$};
\draw[shift={(0,-0.4)}] (2,0) node{\tiny $E_2$};
\draw[shift={(0,-0.4)}] (3,0) node{\tiny $E_3$};
\draw[shift={(0,-0.4)}] (5,0) node{\tiny $E_{l-1}$};
\draw[shift={(0,-0.4)}] (6,0) node{\tiny $E_l$};
\draw[shift={(0,-0.4)}] (8,0) node{\tiny $E_{n-1}$};
\draw[shift={(0,-0.4)}] (9,0) node{\tiny $E_n$};

\draw (E1)--(E2);
\draw (E2)--(E3);
\draw (E3)--(dot);
\draw (dot)--(E_{l-1});
\draw (E_{l-1}.-20)--(E_l.-160);
\draw (E_{l-1})--(E_l);
\draw (E_{l-1}.20)--(E_l.160);
\draw (E_l.-20)--(dot2.-170);
\draw (E_l)--(dot2);
\draw (E_l.20)--(dot2.170);
\draw (dot2)--(E_{n-1});
\draw (dot2.-10)--(E_{n-1}.-160);
\draw (dot2.10)--(E_{n-1}.160);
\draw (E_{n-1}.-20)--(E_n.-160);
\draw (E_{n-1})--(E_n);
\draw (E_{n-1}.20)--(E_n.160);

\end{tikzpicture}
        \]
        and $H^0(E_l, \sO_{E_l})$ is inseparable over $\kappa$.
        \end{itemize}
\end{cor}
\begin{proof}
Proposition \ref{prop:geom red,cusp} handles cusps and twisted cusps, because quadratic extensions are separable in characteristic three.
In the remaining non-elliptic cases, the classification gives $r_i\le4$.
Thus $K_i/\kappa$ can be inseparable only if $r_i=3$.
The graphs containing such a vertex have the displayed form, and the intersections of successive degree-three components identify their constant fields.
The assertion now follows from Proposition \ref{prop:equivalence of geom.red and separability}.
\end{proof}

\begin{eg}\label{eg:geom reducedness of exc div2}
    In Setting \ref{notation:normal surface sings}, suppose that $X=\Spec R$ is klt.
    \begin{enumerate}[label=\textup{(\arabic*)}]
    \item Assume that some exceptional curve $E_i$ is not geometrically reduced over $\kappa$.
    Then by Proposition \ref{prop:equivalence of geom.red and separability} and the classification \cite[Figure 1]{Sato23}, one of the following holds:
    \begin{enumerate}[label=\textup{(\alph*)}]
    \item $p=3$, there is an exceptional curve $E_j$ such that $H^0(E_j, \sO_{E_j})/\kappa$ is inseparable of degree three, and the graph of $\pi$ is the triple folding of $(2,2,2)$. 
    
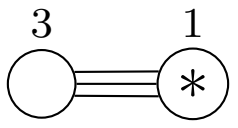
\[
    \begin{tikzpicture}
    \node[draw, shape=circle, inner sep=4.5pt] (B1) at (0,-1){};
    \draw[shift={(0,0.4)}] (0,-1) node{\tiny $3$};
    \node[draw, shape=circle, inner sep=1.8pt] (B2) at (1,-1){$*$};
    \draw[shift={(0,0.4)}] (1,-1) node{\tiny $1$};

    \draw (B1.20)--(B2.160);
    \draw (B1)--(B2);
    \draw (B1.-20)--(B2.200);
    \end{tikzpicture}
    \]
        
    \item $p=2$, there is an exceptional curve $E_j$ such that $H^0(E_j, \sO_{E_j})/\kappa$ is inseparable of degree two, and the graph of $\pi$ is either a twisted chain or a twisted star-shaped graph other than the triple folding of $(2,2,2)$.
    
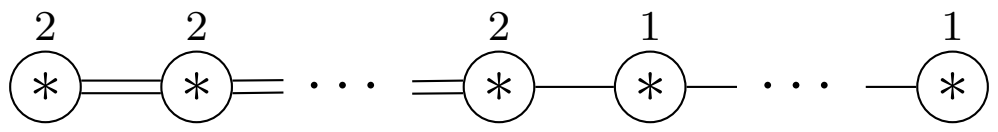
\[
\begin{tikzpicture}
\node[draw, shape=circle, inner sep=1.8pt] (A1) at (0,0){$*$};
\draw[shift={(0,0.4)}] (0,0) node{\tiny $2$};
\node[draw, shape=circle, inner sep=1.8pt] (A2) at (1,0){$*$};
\draw[shift={(0,0.4)}] (1,0) node{\tiny $2$};
\draw node (ADot) at (2,0){$\cdots$};
\node[draw, shape=circle, inner sep=1.8pt] (A3) at (3,0){$*$};
\draw[shift={(0,0.4)}] (3,0) node{\tiny $2$};
\node[draw, shape=circle, inner sep=1.8pt] (A4) at (4,0){$*$};
\draw[shift={(0,0.4)}] (4,0) node{\tiny $1$};
\draw node (ADot2) at (5,0){$\cdots$};
\node[draw, shape=circle, inner sep=1.8pt] (A5) at (6,0){$*$};
\draw[shift={(0,0.4)}] (6,0) node{\tiny $1$};

\draw (A1.10)--(A2.170);
\draw (A1.-10)--(A2.190);
\draw (A2.10)--(ADot.174);
\draw (A2.-10)--(ADot.185);
\draw (ADot.5)--(A3.170);
\draw (ADot.-6)--(A3.-170);
\draw (A3)--(A4);
\draw (A4)--(ADot2);
\draw (ADot2)--(A5);
\end{tikzpicture}
\]
    \item $p=2$ and the graph of $\pi$ is a chain of length one, that is, $\Exc(\pi)$ is an irreducible curve of genus $0$.
    
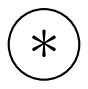
\[
    \begin{tikzpicture}
        \node[draw, shape=circle, inner sep=1.8pt] (A1) at (0,0){$*$};
    \end{tikzpicture}
    \]
    \end{enumerate}
    
    \item Conversely, if either \textup{(a)} or \textup{(b)} holds, then the curve $E_j$ with the specified inseparable constant field is not geometrically reduced over $\kappa$ (Lemma \ref{lem:geom.red to separable}).
\item 
We now give examples of (a), (b), and (c), respectively. In each case, take the local ring at the origin of the displayed affine surface.
In all three examples below, the point blowups are crepant, since their centers are multiplicity-two hypersurface singularities in regular threefolds. 
Thus $K_Y$ is numerically trivial over $X$, the displayed resolutions below are minimal, and the singularities are canonical, in particular klt.

\begin{enumerate}[label=\textup{(\roman*)}]
\item Let $k$ be an $F$-finite field of characteristic $3$ and $a \in k$ be an element with $ \sqrt[3]{a} \not\in k$.
Consider 
\[
X \coloneqq \Spec k[x,y,z]/(z^2+x^3-ay^3+y^7).
\] 
Then the minimal resolution of $X$ is obtained by two successive point blowups, and its exceptional curves are $\PP^1_k$ and $\PP^1_{K}$, where we write $K \coloneqq k(\sqrt[3]{a})$.
Indeed, on the $y$-chart of the first blowup, the strict transform has equation
\[
z_1^2+y(x_1^3-a)+y^5=0.
\]
Its only non-regular point has ideal $(y,z_1,x_1^3-a)$ and residue field $K$; blowing up this point resolves the singularity.
The dual graph at the origin is
    
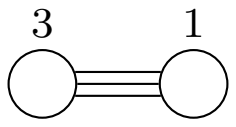
\[
    \begin{tikzpicture}
    \node[draw, shape=circle, inner sep=4.5pt] (B1) at (0,-1){};
    \draw[shift={(0,0.4)}] (0,-1) node{\tiny $3$};
    \node[draw, shape=circle, inner sep=4.5pt] (B2) at (1,-1){};
    \draw[shift={(0,0.4)}] (1,-1) node{\tiny $1$};

    \draw (B1.20)--(B2.160);
    \draw (B1)--(B2);
    \draw (B1.-20)--(B2.200);
    \end{tikzpicture}
    \]
and hence this is an example of (a).
\item Let $k$ be an $F$-finite field of characteristic $2$ and $a, b \in k$ be elements such that $C \coloneqq \Proj k[x,y,z]/(x^2+ay^2+bz^2)$ has no $k$-rational point.
We note that $C$ is a regular projective curve of genus $0$ with $H^0(C, \sO_C) = k$, but not geometrically reduced over $k$.
The minimal resolution of
\[
X \coloneqq \Spec k[x,y,z]/(x^2+ay^2+bz^4).
\] 
is obtained by two successive point blowups, and its exceptional curves are $C$ and $\PP^1_{K}$, where we write $K \coloneqq k(\sqrt{a})$.
On the $z$-chart of the first blowup, the strict transform has equation $x_1^2+ay_1^2+bz^2=0$, the affine cone over $C$.
Blowing up its vertex gives the second exceptional curve $C$ and makes the surface regular.
The dual graph at the origin is
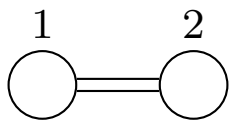
\[
\begin{tikzpicture}
\node[draw, shape=circle, inner sep=4.5pt] (B1) at (0,0){};
\draw[shift={(0,0.4)}] (0,0) node{\tiny $1$};
\node[draw, shape=circle, inner sep=4.5pt] (B2) at (1,0){};
\draw[shift={(0,0.4)}] (1,0) node{\tiny $2$};

\draw (B1.10)--(B2.170);
\draw (B1.-10)--(B2.190);
\end{tikzpicture}
\]
and hence this is an example of (b). Notice that $C$ is not geometrically reduced even though $H^0(C,\sO_C)=k$; the inseparable constant field occurs on the adjacent component $\PP^1_K$.

\item Let $k,a,b$, and $C$ be as in (ii).
The minimal resolution of
\[
X \coloneqq \Spec k[x,y,z]/(x^2+ay^2+bz^2).
\] 
is obtained by blowing up the origin, and its exceptional curve is $C$.
The resulting surface is the total space of $\sO_C(-1)$ and is regular because $C$ is regular.
This gives (c).
\end{enumerate}
\end{enumerate}
\end{eg}

\subsection{Case of \texorpdfstring{$F$}{F}-pure singularities}

Example \ref{eg:geom reducedness of exc div2} shows that, in characteristics $2$ and $3$, exceptional curves on the minimal resolution of an lc surface singularity need not be geometrically reduced.
We now show that $F$-purity rules out this behavior outside the RDP case, and that strong $F$-regularity rules it out altogether.

We first record an observation about Frobenius splittings on a minimal resolution.
In Setting \ref{notation:normal surface sings}, the divisor $K_Y$ is $\pi$-nef, so the negativity lemma gives
\[
\Delta_Y\coloneqq\pi^*K_X-K_Y\geq 0.
\]
Consequently, for every integer $m>0$, restriction identifies
\begin{equation}\label{eq:anticanonical-sections-minimal-resolution}
H^0(Y,\sO_Y(-mK_Y))=H^0(X,\sO_X(-mK_X)).
\end{equation}
These identifications commute with the Frobenius trace maps, since the maps agree away from the exceptional locus.
    \begin{equation}\label{eq:compatible trace}
        \begin{tikzcd}
        H^0(Y, F_*\sO_Y((1-p)K_Y)) \arrow[r,"\cong"] \arrow[d,"\mathrm{Tr}_Y"] & H^0(X , F_*\sO_X((1-p)K_X)) \arrow[d,"\mathrm{Tr}_X"] \\
        H^0(Y,\sO_Y) \arrow[r,"\cong"] & H^0(X,\sO_X)
    \end{tikzcd}
    \end{equation}

\begin{lem}\label{lem:F-split to geom.reduced}
    In the notation of Setting \ref{notation:normal surface sings}, let $E_i$ and $E_j$ be distinct exceptional curves and let $P\in E_i\cap E_j$ be a closed point.
    Suppose that
    \begin{enumerate}
        \item[\textup{(1)}] $E_i$ is geometrically irreducible and has genus zero;
        \item[\textup{(2)}] $K_i=H^0(E_i,\sO_{E_i})$ is separable over $\kappa$; and
        \item[\textup{(3)}] $(E_i,aP)$ is globally sharply $F$-split for some rational number $a>1/p$.
    \end{enumerate}
    Then $E_i$ and $E_j$ are geometrically reduced over $\kappa$, and $K_j=H^0(E_j,\sO_{E_j})$ is separable over $\kappa$.
\end{lem}
\begin{proof}
    By \cite[Corollary 2.5]{GLP$^+$15}, assumptions (2) and (3) imply that $E_i$ is geometrically reduced over $\kappa$.
    Together with (1), this shows that $E_i$ is smooth over $\kappa$ \cite[Lemma 3.1]{Sato23}.
    Theorem \ref{thm:F-split to geom F-split} therefore implies that $(E_i,aP)$ is geometrically sharply $F$-split.
    Let $q$ be the inseparable degree of $\kappa(P)/\kappa$.
    After base change to an algebraic closure, the pullback of $P$ has multiplicity $q$ at each point in its support.
    Sharp $F$-purity forces every boundary coefficient to be at most one, so $aq\leq 1$.
    Since $q$ is a power of $p$ and $a>1/p$, we obtain $q=1$.
    Thus $\kappa(P)/\kappa$ is separable.

    The scheme-theoretic intersection $E_i\cap E_j$ is reduced: this follows from the intersection numbers in the classification of dual graphs, with the two-component cusp case covered by Proposition \ref{prop:geom red,cusp}.
    Hence $P$ is an effective Cartier divisor on $E_j$.
    Lemma \ref{lem:geom.red} now shows that $E_j$ is geometrically reduced over $\kappa$, and Lemma \ref{lem:geom.red to separable} gives the separability of $K_j/\kappa$.
\end{proof}

\begin{prop}\label{prop:geom red,p=3,SFP}
    In the notation of Setting \ref{notation:normal surface sings}, suppose that $p=3$ and $X$ is $F$-pure.
    Then every exceptional curve $E_i$ is geometrically reduced over $\kappa$, with the possible exception of an RDP with the following dual graph:
\[
\begin{tikzpicture}
\node[draw, shape=circle, inner sep=4.5pt] (E1) at (3,0){};
\node[draw, shape=circle, inner sep=4.5pt] (E2) at (4,0){};
\draw[shift={(0,0.4)}] (3,0) node{\tiny $3$};
\draw[shift={(0,0.4)}] (4,0) node{\tiny $1$};
\draw[shift={(0,-0.4)}] (3,0) node{\tiny $E_1$};
\draw[shift={(0,-0.4)}] (4,0) node{\tiny $E_2$};

\draw (E1)--(E2);
\draw (E1.10)--(E2.170);
\draw (E1.-10)--(E2.-170);
\end{tikzpicture}
\]
\end{prop}
\begin{proof}
    By \cite[Lemma 4.3]{Sato23}, every finite separable extension $l/\kappa$ is realized as the residue field of a finite \'etale local extension $R'/R$.
    Since $F$-purity is preserved under finite \'etale morphisms \cite[Theorem 4.8]{Hara-Watanabe}, we may make such a base change and assume that every exceptional curve is geometrically irreducible over $\kappa$.
    Geometric reducedness descends along the resulting field extension.

    Suppose first that $X$ is simple elliptic, and let $E$ be its exceptional curve.
    Since $\pi^*K_X=K_Y+E$, restriction gives an isomorphism
    \[
    H^0(Y,\sO_Y((1-p)(K_Y+E)))\cong H^0(X,\sO_X((1-p)K_X)).
    \]
    Compatibility of the trace maps \eqref{eq:compatible trace} and $F$-purity of $X$ show that $(Y,E)$ is globally sharply $F$-split.
    By Remark \ref{rem:non-normal adjunction}, $E$ is globally sharply $F$-split.
    Since $H^0(E,\sO_E)=\kappa$ (\cite[Corollary 10.10]{Kol13}), it follows from \cite[Corollary 2.5]{GLP$^+$15} that $E$ is geometrically reduced over $\kappa$.

    We may now assume that $X$ is not simple elliptic.
    By Corollary \ref{cor:geom red,p=3}, it remains to consider the following graph, with $K_l/\kappa$ inseparable:
    
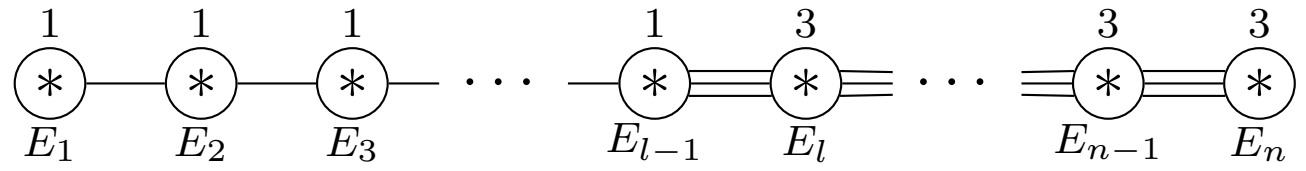
\[
\begin{tikzpicture}
\node[draw, shape=circle, inner sep=1.8pt] (E1) at (1,0){$*$};
\node[draw, shape=circle, inner sep=1.8pt] (E2) at (2,0){$*$};
\node[draw, shape=circle, inner sep=1.8pt] (E3) at (3,0){$*$};
\draw node (dot) at (4,0){$\cdots$};
\node[draw, shape=circle, inner sep=1.8pt] (E_{l-1}) at (5,0){$*$};
\node[draw, shape=circle, inner sep=1.8pt] (E_l) at (6,0){$*$};
\draw node (dot2) at (7,0){$\cdots$};
\node[draw, shape=circle, inner sep=1.8pt] (E_{n-1}) at (8,0){$*$};
\node[draw, shape=circle, inner sep=1.8pt] (E_n) at (9,0){$*$};

\draw[shift={(0,0.4)}] (1,0) node{\tiny $1$};
\draw[shift={(0,0.4)}] (2,0) node{\tiny $1$};
\draw[shift={(0,0.4)}] (3,0) node{\tiny $1$};
\draw[shift={(0,0.4)}] (4,0) node{\tiny ${}$};
\draw[shift={(0,0.4)}] (5,0) node{\tiny $1$};
\draw[shift={(0,0.4)}] (6,0) node{\tiny $3$};
\draw[shift={(0,0.4)}] (8,0) node{\tiny $3$};
\draw[shift={(0,0.4)}] (9,0) node{\tiny $3$};

\draw[shift={(0,-0.4)}] (1,0) node{\tiny $E_1$};
\draw[shift={(0,-0.4)}] (2,0) node{\tiny $E_2$};
\draw[shift={(0,-0.4)}] (3,0) node{\tiny $E_3$};
\draw[shift={(0,-0.4)}] (5,0) node{\tiny $E_{l-1}$};
\draw[shift={(0,-0.4)}] (6,0) node{\tiny $E_l$};
\draw[shift={(0,-0.4)}] (8,0) node{\tiny $E_{n-1}$};
\draw[shift={(0,-0.4)}] (9,0) node{\tiny $E_n$};

\draw (E1)--(E2);
\draw (E2)--(E3);
\draw (E3)--(dot);
\draw (dot)--(E_{l-1});
\draw (E_{l-1}.-20)--(E_l.-160);
\draw (E_{l-1})--(E_l);
\draw (E_{l-1}.20)--(E_l.160);
\draw (E_l.-20)--(dot2.-170);
\draw (E_l)--(dot2);
\draw (E_l.20)--(dot2.170);
\draw (dot2)--(E_{n-1});
\draw (dot2.-10)--(E_{n-1}.-160);
\draw (dot2.10)--(E_{n-1}.160);
\draw (E_{n-1}.-20)--(E_n.-160);
\draw (E_{n-1})--(E_n);
\draw (E_{n-1}.20)--(E_n.160);

\end{tikzpicture}
    \]

    For $l\leq j<n$, the scheme $E_j\cap E_{j+1}$ is a scheme over both $K_j$ and $K_{j+1}$, and its coordinate ring has dimension $3$ over $\kappa$.
    Since $[K_j:\kappa]=[K_{j+1}:\kappa]=3$, we obtain
    \[
    E_j\cap E_{j+1}\cong\Spec K_j\cong\Spec K_{j+1}.
    \]
    Thus $K_l\cong\cdots\cong K_n$ as $\kappa$-algebras, and it suffices to prove that $K_l/\kappa$ is separable.
    Set $E\coloneqq\sum_{i=1}^n E_i$.

    \begin{cl}
        The pair $(Y,\frac12(E-E_{l-1})+E_{l-1})$ is globally sharply $F$-split.
    \end{cl}
    \begin{proof}[Proof of the claim]
        Since $X$ is not an RDP, there is an exceptional curve $E_i$ with $E_i^2<-2r_i$.
        Adjunction gives
        \[
        2K_Y\cdot E_i=-4r_i-2E_i^2>0.
        \]
        Thus $\sO_{E_i}(-2K_Y)$ has negative degree, and the restriction sequence gives
        \[
        H^0(Y,\sO_Y(-2K_Y-E_i))=H^0(Y,\sO_Y(-2K_Y)).
        \]
        If $Z$ is a reduced exceptional divisor already subtracted in this way and $E_j\not\subseteq\Supp Z$ meets $Z$, then
        \[
        (2K_Y+Z)\cdot E_j>0,
        \]
        because $K_Y$ is $\pi$-nef.
        Subtracting adjacent components successively and using connectedness, we obtain
        \[
        H^0(Y,\sO_Y(-2K_Y-E))=H^0(Y,\sO_Y(-2K_Y)).
        \]
        Finally,
        \[
        \begin{aligned}
        (2K_Y+E)\cdot E_{l-1}
        &=2(K_Y+E_{l-1})\cdot E_{l-1}-E_{l-1}^2
          +(E-E_{l-1})\cdot E_{l-1}\\
        &\geq -4+2+3=1.
        \end{aligned}
        \]
        Another restriction sequence therefore gives
        \[
        H^0(Y,\sO_Y(-2K_Y-E-E_{l-1}))
        =H^0(Y,\sO_Y(-2K_Y)).
        \]
        By \eqref{eq:anticanonical-sections-minimal-resolution}, \eqref{eq:compatible trace} and $F$-purity of $X$, the Frobenius trace on these sections is surjective.
        Since $p-1=2$, this is exactly the asserted global sharp $F$-splitting.
    \end{proof}

    The genus-zero curve $E_{l-1}$ is smooth over $\kappa$, since $K_{l-1}=\kappa$, $p=3$, and it is geometrically irreducible \cite[Lemmas 3.1 and 3.3]{Sato23}.
    Adjunction (Lemma \ref{lem:adjunction}) and removal of the other boundary components show that
    $(E_{l-1},\frac12 E_l|_{E_{l-1}})$ is globally sharply $F$-split.
    Here $K_{l-1}=\kappa$ and $E_{l-1}\cap E_l\cong\Spec K_l$.
    Applying Lemma \ref{lem:F-split to geom.reduced}, with $a=1/2>1/3$, shows that $K_l/\kappa$ is separable, as required.
\end{proof}

We next consider characteristic $2$.

\begin{prop}\label{prop:geom red,p=2,SFP}
    In the notation of Setting \ref{notation:normal surface sings}, suppose that $p=2$ and $X$ is $F$-pure but not an RDP.
    Then every exceptional curve is geometrically reduced over $\kappa$.
\end{prop}
\begin{proof}
    As in the proof of Proposition \ref{prop:geom red,p=3,SFP}, we may assume that every exceptional curve is geometrically irreducible over $\kappa$.
    Set $E\coloneqq\sum_i E_i$.
    We first show that $(Y,E)$ is globally sharply $F$-split.
    If $X$ is not rational, then $\pi^*K_X=K_Y+E$ \cite[Proposition 2.28]{Kol13}, so
    \[
    H^0(Y,\sO_Y(-K_Y-E))=H^0(X,\sO_X(-K_X)).
    \]
    If $X$ is rational, all $E_i$ have genus zero, and the assumption that $X$ is not an RDP gives a component with $K_Y\cdot E_i>0$.
    The same successive restriction argument used in the preceding proposition, now with $K_Y$ in place of $2K_Y$, gives
    \[
    H^0(Y,\sO_Y(-K_Y-E))=H^0(Y,\sO_Y(-K_Y))
    =H^0(X,\sO_X(-K_X)).
    \]
    In either case, compatibility with Frobenius trace and $F$-purity of $X$ prove the claim.
    In particular, each component $E_i$ is Frobenius split by Remark \ref{rem:non-normal adjunction}.
    If $E$ has only one irreducible component, then noting that $H^0(E,\sO_E)=\kappa$ (\cite[Corollary 10.10]{Kol13}), the assertion follows from \cite[Corollary 2.5]{GLP$^+$15}.
    We may therefore assume that $E$ has at least two components, all of genus zero.
    We also note that \cite[Corollary 2.5]{GLP$^+$15} and \cite[Lemma 3.1]{Sato23} imply that each $E_i$ is smooth over $K_i$.

    Take exceptional curves $E_i$ and $E_j$ such that $E_i \cap E_j \neq \emptyset$.
    Noting that $(E_i, E_j|_{E_i})$ is globally sharply $F$-split, it follows from the proof of Lemma \ref{lem:F-split to geom.reduced} that $\kappa(P)$ is separable over $H^0(E_i,\sO_{E_i})$ for every point $P \in E_i \cap E_j$.
    Combining this with Proposition \ref{prop:geom red,cusp}, we conclude the assertion in the case of a cusp or a twisted cusp.

    In the remaining cases, the classification provides a component $E_1$ with $K_1=\kappa$ (\cite[Corollary 10.10]{Kol13}).
    Lemma \ref{lem:F-split to geom.reduced}, with $a=1$, gives geometric reducedness of each adjacent component and separability of its constant field over $\kappa$.
    Connectedness of the exceptional divisor completes the induction.
\end{proof}

\begin{cor}\label{cor:SFP to Geom red}
    In the notation of Setting \ref{notation:normal surface sings}, suppose that $X$ is $F$-pure but not an RDP.
    Then every $E_i$ is geometrically reduced over $\kappa$.
\end{cor}
\begin{proof}
    Combine Corollary \ref{cor:geom red,p>3} with Propositions \ref{prop:geom red,p=3,SFP} and \ref{prop:geom red,p=2,SFP}.
\end{proof}

For strongly $F$-regular singularities, the RDP cases can also be included.

\begin{thm}\label{thm:SFR to Geom.red}
    In the notation of Setting \ref{notation:normal surface sings}, suppose that $X$ is strongly $F$-regular.
    Then every $E_i$ is geometrically reduced over $\kappa$.
\end{thm}
\begin{proof}
    Strong $F$-regularity is preserved under finite \'etale morphisms \cite[Theorem 4.8]{Hara-Watanabe}.
    Thus we may again assume that every exceptional curve is geometrically irreducible over $\kappa$.
    By Corollary \ref{cor:geom red,p>3}, only $p=2,3$ remain.

    Suppose first that $p=3$.
    By Proposition \ref{prop:geom red,p=3,SFP}, it remains to consider an RDP with the following dual graph:
    
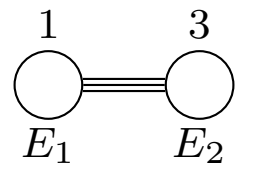
\[
\begin{tikzpicture}
\node[draw, shape=circle, inner sep=4.5pt] (E1) at (3,0){};
\node[draw, shape=circle, inner sep=4.5pt] (E2) at (4,0){};
\draw[shift={(0,0.4)}] (3,0) node{\tiny $1$};
\draw[shift={(0,0.4)}] (4,0) node{\tiny $3$};
\draw[shift={(0,-0.4)}] (3,0) node{\tiny $E_1$};
\draw[shift={(0,-0.4)}] (4,0) node{\tiny $E_2$};

\draw (E1)--(E2);
\draw (E1.10)--(E2.170);
\draw (E1.-10)--(E2.-170);
\end{tikzpicture}
    \]
    Since $X$ is Gorenstein, Lemma \ref{lem:key} provides a nonzero $x\in\m$ such that $(X,\frac12 D)$ is sharply $F$-pure, where $D\coloneqq\operatorname{div}(x)$.
    The equality $\pi^*K_X=K_Y$ and compatibility of trace maps show that $(Y,\frac12\pi^*D)$ is globally sharply $F$-split.
    Write
    \[
    \pi^*D=D'+b_1E_1+b_2E_2,
    \]
    where $D'$ is the strict transform of $D$.
    The $b_i$ are positive integers, and sharp $F$-purity gives $b_i\leq2$.
    Intersecting with $E_1$ yields
    \[
    0=D'\cdot E_1-2b_1+3b_2.
    \]
    Since $D'\cdot E_1\geq0$ and $b_i\in\{1,2\}$, we obtain $(b_1,b_2)=(2,1)$.
    Thus $(Y,E_1+\frac12 E_2)$ is globally sharply $F$-split.
    Adjunction and Lemma \ref{lem:F-split to geom.reduced} prove geometric reducedness of both components.

    Now suppose that $p=2$.
    By Proposition \ref{prop:geom red,p=2,SFP}, we may assume that $X$ is an RDP.
    Lemma \ref{lem:key} provides a nonzero $x\in\m$ such that $(Y,\pi^*D)$ is globally sharply $F$-split, where $D\coloneqq\operatorname{div}(x)$.
    Since $\pi^*D\geq E\coloneqq\sum_i E_i$, the pair $(Y,E)$ is also globally sharply $F$-split.
    Then, as in the proof of Proposition \ref{prop:geom red,p=2,SFP}, the repeated use of Lemma \ref{lem:F-split to geom.reduced} shows the assertion.
\end{proof}

\begin{cor}\label{cor:minres=logres}
    In the notation of Setting \ref{notation:normal surface sings}, suppose that $X$ is strongly $F$-regular.
    Then there is a finite \'etale cover $X'=\Spec R'\to X$, with $R'$ local, such that the exceptional curves $E_i'$ of the minimal resolution
    \[
    \pi'\coloneqq\pi\times_X\id_{X'}\colon Y\times_X X'\to X'
    \]
    satisfy the following properties, where $K$ is the residue field of $R'$:
    \begin{enumerate}[label=\textup{(\roman*)}]
        \item $E_i'\cong\PP^1_K$ for every $i$;
        \item for $i\neq j$, the scheme-theoretic intersection $E_i'\cap E_j'$ is either empty or isomorphic to $\Spec K$, and no three components meet; in particular, the exceptional divisor is a simple normal crossing curve over $K$;
        \item the dual graph is non-twisted.
    \end{enumerate}
\end{cor}
\begin{proof}
    After a finite \'etale local base change, we may assume that every exceptional curve is geometrically irreducible over $K$.
    By Theorem \ref{thm:SFR to Geom.red}, every $E_i'$ is then geometrically integral, so $H^0(E_i',\sO_{E_i'})=K$.
    This proves (iii).
    Each $E_i'$ has genus zero, and Lemma \ref{lem:genus zero curve}(2) allows a further finite separable extension of $K$ over which every component is isomorphic to $\PP^1$.
    Realizing this extension by a finite \'etale local cover proves (i).
    Finally, the classification of klt dual graphs \cite[Figure 1]{Sato23} gives (ii).
\end{proof}

\section{Classification over imperfect residue fields}\label{section:Appendix}
In this section, we prove Theorems \ref{Introthm:classification of SFR} and \ref{Introthm:classification of SFP}.

\subsection{Strongly \texorpdfstring{$F$}{F}-regular rational double points}

We use Setting \ref{notation:normal surface sings} throughout this subsection.
For a strongly $F$-regular rational double point, we establish the bounds of the characteristic $p$ associated with the non-twisted star-shaped graphs.
Over an algebraically closed field, Hara obtains these bounds from explicit equations and Fedder's criterion.
We instead use Lemma \ref{lem:key} and the intersection matrix of the exceptional curves, so the argument also applies over an imperfect residue field.

We first record the common setup for the three proofs below.
Since an RDP is Gorenstein, Lemma \ref{lem:key} provides a nonzero element $f\in\m$ such that, for $D=\operatorname{div}(f)$, the pair $(X,D/(p-1))$ is sharply $F$-pure.
The minimal resolution is crepant, so $(Y,\pi^*D/(p-1))$ is globally sharply $F$-split.
Write
\[
\pi^*D=D'+\sum_i a_iE_i,
\]
where $D'$ is the strict transform and $a_i\in\mathbb Z_{>0}$.
Sharp $F$-purity at the generic point of $E_i$ gives $a_i\le p-1$.
Moreover, because $\pi^*D\cdot E_i=0$ and $D'\cdot E_i\ge0$, the exceptional part is anti-nef.
For a non-twisted RDP graph, this means
\begin{equation}\label{eq:RDP-antinef}
2a_i\ge\sum_{E_j\text{ adjacent to }E_i}a_j.
\end{equation}

\begin{lem}\label{lem:RDP1}
Suppose that $X$ is a strongly $F$-regular RDP and that the dual graph of its minimal resolution is star-shaped.
Then $p>2$.
\end{lem}
\begin{proof}
Let $a$ be the coefficient of the central curve in $\pi^*D$, and let $b_1,b_2,b_3$ be the coefficients of the three adjacent curves.
By \eqref{eq:RDP-antinef},
\[
2a\ge b_1+b_2+b_3\ge3.
\]
Thus $a\ge2$, and $p-1\ge a$ gives $p>2$.
\end{proof}

\begin{lem}\label{lem:RDP2}
Suppose that $X$ is a strongly $F$-regular RDP and that the dual graph of its minimal resolution is star-shaped of type $(2,3,3)$ or $(2,3,4)$.
Then $p>3$.
\end{lem}
\begin{proof}
In either case, 
two arms have length at least two.
Let $a$ be the coefficient of the central curve, let $b_1,b_2$ be the coefficients of its neighbors on these two arms, and let $c_1,c_2$ be the coefficients of the next curves on the same arms.
Let $b_3$ be the coefficient of the neighbor on the third arm.
Then
\[
2a\ge b_1+b_2+b_3,\qquad
2b_1\ge a+c_1,\qquad
2b_2\ge a+c_2.
\]
Lemma \ref{lem:RDP1} and its proof give $a\ge2$.
If $a=2$, then $c_1,c_2\ge1$ imply $b_1,b_2\ge2$, contradicting $2a\ge b_1+b_2+b_3$.
Consequently $a\ge3$, so $p-1\ge3$ and the primality of $p$ gives $p\ge5$.
\end{proof}

\begin{lem}\label{lem:RDP3}
Suppose that $X$ is a strongly $F$-regular RDP and that the dual graph of its minimal resolution is star-shaped of type $(2,3,5)$.
Then $p>5$.
\end{lem}
\begin{proof}
The three arms of the $E_8$ graph have lengths $1$, $2$, and $4$.
Let $a$ be the coefficient of the central curve.
For an arm of length $\ell$, write its coefficients, starting next to the central curve, as $b_1,\ldots,b_\ell$, and put $b_0=a$ and $b_{\ell+1}=0$.
Multiplying the inequalities
\[
2b_i-b_{i-1}-b_{i+1}\ge0\qquad(1\le i\le\ell)
\]
by $\ell+1-i$ and summing gives
\[
(\ell+1)b_1\ge\ell a.
\]
Applying this to the three arms and then using \eqref{eq:RDP-antinef} at the central curve yields
\[
2a\ge
\left\lceil\frac a2\right\rceil+
\left\lceil\frac{2a}{3}\right\rceil+
\left\lceil\frac{4a}{5}\right\rceil.
\]
For $a=1,2,3,4,5$, the right-hand side equals $2a+1$.
Hence $a\ge6$, and $p-1\ge a$ implies $p\ge7$.
\end{proof}

\subsection{Proofs of the main theorems}

\begin{proof}[Proof of Theorem \ref{Introthm:classification of SFR}]
The regular case is immediate, so assume that $S$ is singular.
Strong $F$-regularity is preserved under finite \'etale covers and descends under faithfully flat maps
\cite[Corollary 2.4]{Hochster-Huneke89}, \cite[Theorem 4.8]{Hara-Watanabe}.
If $S$ is strongly $F$-regular, then it is klt by
\cite[Theorem 3.3]{Hara-Watanabe} and \cite[Theorem 4.3]{SS10},
and Condition (A) follows from Theorem \ref{thm:SFR to Geom.red}.
In the other direction, Condition (B) implies that $S$ is klt.
Thus we may assume that $S$ is klt and Condition (A) holds in both directions of the argument.

As in the proof of Corollary \ref{cor:minres=logres}, after a finite \'etale local cover we may arrange that all exceptional curves are isomorphic to $\PP^1_\kappa$ and all their intersection points are $\kappa$-rational.
Lemma \ref{lem:geom.irr} preserves the relevant graph types under this replacement.

Suppose first that Conditions (A) and (B) hold.
Proposition \ref{prop:log Fano curve} supplies the curve-theoretic input in Hara's proof, and the argument of \cite[Section 4.2]{Hara(two-dim)} gives strong $F$-regularity.
We spell out why the vanishing used in that argument remains valid over our residue field.
Using Hara's notation $X,D,E_0$ from \cite[Section 3.8]{Hara(two-dim)} and writing $q=p^e$, the required vanishing is
\[
H^1\bigl(X,\omega_X^{1-q}(-\lfloor qD\rfloor+nE_0)\bigr)
=H^1\bigl(X,\sO_X(K_X+\lceil nE_0-q(K_X+D)\rceil)\bigr)=0.
\]
Put $N=nE_0-q(K_X+D)$ and $\Delta_X=\lceil N\rceil-N$.
For $q$ sufficiently large, \cite[Lemma 3.9]{Hara(two-dim)} gives $N\cdot E_i\ge0$ for every exceptional curve $E_i$.
Moreover, $\Delta_X$ is supported on the simple normal crossing exceptional locus and has all coefficients in $[0,1)$.
Thus $\lceil N\rceil-\Delta_X=N$ is relatively nef, and the required vanishing follows from \cite[Theorem 10.4]{Kol13}.

Conversely, suppose that $S$ is strongly $F$-regular.
The classification of klt dual graphs reduces Condition (B) to the bounds of characteristic for the three star-shaped types.
For non-RDPs, the argument of \cite[Section 4.4]{Hara(two-dim)} uses the intersection matrix and the criterion of global $F$-regularity on the central curve; Proposition \ref{prop:log Fano curve} supplies the latter over $\kappa$.
For RDPs, the required bounds follow from Lemmas \ref{lem:RDP1}, \ref{lem:RDP2}, and \ref{lem:RDP3}.
\end{proof}

\begin{proof}[Proof of Theorem \ref{Introthm:classification of SFP}]
$F$-purity is preserved under finite \'etale covers and descends under faithfully flat maps
\cite[Proposition 5.13]{HR}, \cite[Theorem 4.8]{Hara-Watanabe}.

We first prove (II).
Case (1) is immediate.
In cases (3-i), (3-ii), (4), and (5) with $p>2$, the classification of dual graphs, together with Corollaries \ref{cor:geom red,p>3} and \ref{cor:geom red,p=3}, shows that every exceptional curve is geometrically reduced over $\kappa$.
For the non-cusp cases, a finite \'etale local cover, as in the proof of Corollary \ref{cor:minres=logres}, reduces to exceptional curves isomorphic to $\PP^1_\kappa$ with $\kappa$-rational intersection points.
Proposition \ref{prop:log Fano curve} then supplies the splitting criterion needed in \cite[Section 4.2]{Hara(two-dim)}.
The vanishing used on p.~48 of that paper is valid here as well: with Hara's notation,
\[
H^1\bigl(X,\omega_X^{1-p}(-\lfloor p\Delta\rfloor)\bigr)
=H^1\bigl(X,\sO_X(K_X+\lceil-p(K_X+\Delta)\rceil)\bigr)=0
\]
by \cite[Theorem 10.4]{Kol13}: here $-p(K_X+\Delta)$ is relatively numerically trivial, and its rounding error is supported on the simple normal crossing exceptional divisor with coefficients in $[0,1)$.

For a cusp, the argument of \cite[Proposition 4.7]{Hara(two-dim)} reduces the assertion to global $F$-splitting of the reduced exceptional curve $\Delta_0 = \bigcup_i E_i$.
Noting that we have $K_{\Delta_0} = (K_T+ \Delta)|_{\Delta_0} \sim 0$, by the similar argument as in the proof of \cite[Claim 4.8.1]{Hara(two-dim)}, we may reduce it to the globally sharp $F$-splitting of the pair $(\PP^1_{\kappa}, [0]+[\infty])$, which is true in any characteristic.

For (I), assume that $S$ is $F$-pure.
Then $S$ is lc by \cite[Theorem 3.3]{Hara-Watanabe} and \cite[Theorem 4.3]{SS10}.
The assertions are immediate if $S$ is strongly $F$-regular, an RDP, or simple elliptic.
In the remaining cases, Corollary \ref{cor:SFP to Geom red} gives geometric reducedness of the exceptional curves.
The classification of lc dual graphs places $S$ among the non-RDPs klt or the graphs listed in (3)--(5).
After a finite \'etale local cover splitting the genus-zero components and their intersection points, the argument of \cite[Sections 4.2 and 4.4]{Hara(two-dim)}, together with Proposition \ref{prop:log Fano curve}, gives the stated characteristic restrictions and excludes non-RDPs klt that are not strongly $F$-regular.
\end{proof}

\begin{rem}\label{rem:Remark on Theorem B}
\
\begin{enumerate}[label=\textup{(\arabic*)}]
\item In case (5) of Theorem \ref{Introthm:classification of SFP} (II), the assumption $p>2$ is used to ensure geometric reducedness of the exceptional curves.
In characteristic $2$, the same proof applies whenever every exceptional curve is geometrically reduced.
Conversely, if an exceptional curve fails to be geometrically reduced, Corollary \ref{cor:SFP to Geom red} shows that $S$ is not $F$-pure.
Proposition \ref{prop:geom red,cusp} describes this obstruction in terms of inseparable residue-field extensions.
Example \ref{eg:non F-pure cusp} gives a non-$F$-pure twisted cusp.

\item Consider cases (3-iii) and (6), which are not included in the converse statement.
Let $\pi\colon T\to S$ be the minimal resolution, and define the effective exceptional $\Q$-divisor $\Delta_T$ by
\[
K_T+\Delta_T=\pi^*K_S.
\]
Since $S$ is lc but not klt, $Z=\lfloor\Delta_T\rfloor$ is nonzero and reduced.
Put $\delta=\{\Delta_T\}|_Z$.
The argument of \cite[Proposition 4.7]{Hara(two-dim)} identifies $F$-purity of $S$ with global $F$-splitting of $(Z,\delta)$.
Example \ref{eg:Theorem B for D4} makes the criterion explicit for $\widetilde D_4$.

In case (6), since we have $\Delta_T=0$, this says that $S$ is $F$-pure exactly when its exceptional curve $E$ is globally $F$-split.
(If $E$ is nonnormal, global $F$-splitting means splitting of the usual map $\sO_Z\to F_*\sO_Z$; no normality assumption is imposed on this map.)
\end{enumerate}
\end{rem}

\begin{eg}\label{eg:non F-pure cusp}
Let $k$ be an $F$-finite field of characteristic $2$, and choose $a\in k\setminus k^2$.
We construct a non-$F$-pure twisted cusp at the origin of
\[
X=\Spec k[x,y,z]/(x^2+y^7+z^7-ay^2z^2).
\]
This hypersurface is integral and regular away from the origin, hence normal.
If $\mathfrak m=(x,y,z)$, then
\[
x^2+y^7+z^7-ay^2z^2\in\mathfrak m^{[2]}.
\]
Fedder's criterion therefore shows that $X$ is not $F$-pure at the origin.

Let
\[
X'=\Proj\left(\bigoplus_{n\ge0}I_nT^n\right)\longrightarrow X,
\qquad I_n=(x^iy^jz^l\mid 2i+j+l\ge n),
\]
be the weighted blow-up with weights $(2,1,1)$.
It is covered by $D_+(yT)$ and $D_+(zT)$, each isomorphic to
\[
V=\Spec k[u,v,w]/\bigl(u^2+v^3(1+w^7)-aw^2\bigr).
\]
The surface $V$ is regular away from the origin.
Blowing up the origins of these two charts gives a regular surface $Y$.
The three exceptional curves of $Y\to X$ are
\[
E_1\cong\PP^1_{k(\sqrt a)},\qquad
E_2\cong E_3\cong\Proj k[U,V,W]/(U^2-aW^2).
\]
Their intersection numbers are
\[
E_1^2=-6,\qquad E_2^2=E_3^2=-2,\qquad
E_1\cdot E_2=E_1\cdot E_3=2,\qquad E_2\cdot E_3=0.
\]
In particular, no exceptional curve is a curve of the first kind, so $Y\to X$ is the minimal resolution.
Its dual graph is the twisted cusp as follows:
    
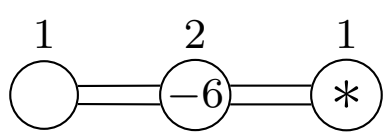
\[
    \begin{tikzpicture}
    \node[draw, shape=circle, inner sep=4.5pt] (B1) at (0,0){};
    \draw[shift={(0,0.4)}] (0,0) node{\tiny $1$};
    \node[draw, shape=circle, inner sep=0.3pt] (B2) at (1,0){\tiny{$-6$}};
    \draw[shift={(0,0.4)}] (1,0) node{\tiny $2$};
    \node[draw, shape=circle, inner sep=1.8pt] (B3) at (2,0){$*$};
    \draw[shift={(0,0.4)}] (2,0) node{\tiny $1$};

    \draw (B1.15)--(B2.165);
    \draw (B1.-15)--(B2.-165);
    \draw (B2.15)--(B3.165);
    \draw (B2.-15)--(B3.-165);
    \end{tikzpicture}
    \]

Finally, the blow-up $Y'\to Y$ at the singular points of $E_2$ and $E_3$ gives a log resolution of $X$. 
Every exceptional coefficient of $K_{Y'/X}$ is therefore $-1$, proving that $X$ is log canonical.

\end{eg}

\begin{eg}\label{eg:Theorem B for D4}
Let $S$ be an lc surface singularity with non-twisted dual graph $\widetilde D_4$:
            
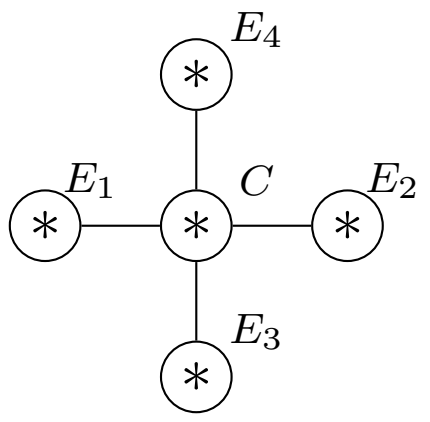
\[
\begin{tikzpicture}
\node[draw, shape=circle, inner sep=1.8pt] (E1 alpha) at (2,0){$*$};
\node[draw, shape=circle, inner sep=1.8pt] (C) at (3,0){$*$};
\node[draw, shape=circle, inner sep=1.8pt] (E1 Beta) at (4,0){$*$};
\node[draw, shape=circle, inner sep=1.8pt] (E1 gamma) at (3,-1){$*$};
\node[draw, shape=circle, inner sep=1.8pt] (E1 delta) at (3,1){$*$};

\draw[shift={(0.3,0.3)}] (2,0) node{\tiny $E_1$};
\draw[shift={(0.4,0.3)}] (3,0) node{\tiny $C$};
\draw[shift={(0.3,0.3)}] (4,0) node{\tiny $E_2$};
\draw[shift={(0.4,0.3)}] (3,-1) node{\tiny $E_3$};
\draw[shift={(0.4,0.3)}] (3,1) node{\tiny $E_4$};

\draw (E1 alpha)--(C);
\draw (C)--(E1 Beta);
\draw (C)--(E1 gamma);
\draw (C)--(E1 delta);
\end{tikzpicture}
        \]

Write $C$ for the central curve and $E_1,\ldots,E_4$ for its four neighbors.
Since $C$ has genus zero and carries a Cartier divisor of degree one, $C\cong\PP^1_\kappa$.
The equalities $C\cdot E_i=1$ show that $C$ and $E_i$ meet transversely at distinct $\kappa$-rational points $P_i$.
The pair in Remark \ref{rem:Remark on Theorem B} (2) is
\[
(Z,\delta)\cong
\left(\PP^1_\kappa,\frac12(P_1+P_2+P_3+P_4)\right).
\]
Choose a coordinate with $P_1=0$, $P_2=\infty$, $P_3=-1$, and $P_4=\lambda$, where $\lambda\in\kappa\setminus\{0,-1\}$.
Proposition \ref{prop:log Fano curve}, equivalently \cite[Theorem 1.2(3)(vi)]{Hara(two-dim)}, then gives the following criterion:
$S$ is $F$-pure if and only if $p=2n+1$ for some integer $n\ge1$ and the coefficient of $x^n$ in
\[
(x+1)^n(x-\lambda)^n
\]
is nonzero.
\end{eg}

The authors learned from Austyn Simpson that the following corollary follows in any dimension from \cite[Theorem 3.6]{Aberbach}.
We give a proof in dimension two using Theorem \ref{Introthm:classification of SFR}.

\begin{cor}
Let $(R,\m,\kappa)$ be a two-dimensional Noetherian $F$-finite normal local domain containing a field $k$ of characteristic $p>0$.
Assume that $X=\Spec R$ is geometrically normal over $k$, that $\kappa/k$ is separable, and that $X$ is strongly $F$-regular.
Then $X$ is geometrically strongly $F$-regular over $k$.
\end{cor}
\begin{proof}
Let $K/k$ be a finite field extension.
After making a finite separable extension and localizing at each closed point, we may assume that $K/k$ is purely inseparable: strong $F$-regularity is preserved under finite \'etale maps \cite[Theorem 4.8]{Hara-Watanabe}.
Since $\kappa/k$ is separable,
\[
L=\kappa\otimes_kK
\]
is a field \cite[Theorem 26.4]{Matsumura}.
Consequently $X_K$ is local with residue field $L$.
If $X$ is regular, the flat local map $R\to R\otimes_kK$ has regular closed fiber $L$, so $X_K$ is regular as well.
We may therefore assume that $X$ is singular.

Take a finite \'etale local cover $X'\to X$ as in Corollary \ref{cor:minres=logres}.
Its residue field is separable over $k$, and $X'$ remains geometrically normal and strongly $F$-regular.
Since strong $F$-regularity descends along the faithfully flat map $X'_K\to X_K$, it suffices to work with $X'$.
We may therefore assume that, for the minimal resolution $\pi\colon Y\to X$, every exceptional curve is $\PP^1_\kappa$ and all intersections are transverse and $\kappa$-rational.

We show that $Y_K$ is regular.
At a point above the exceptional locus, choose an exceptional component through its image.
Its base change is a regular Cartier divisor on $Y_K$, so regularity of the ambient local ring follows from deformation of regularity \stacksproj{00NU}.
Away from the exceptional locus, $Y_K$ is isomorphic to the punctured spectrum of $X_K$.
This is regular because $X_K$ is a normal local surface.
Thus
\[
\pi_K\colon Y_K\longrightarrow X_K
\]
is a resolution.
Its exceptional curves and intersections are obtained from those of $\pi$ by the field extension $\kappa\subset L$.
Their self-intersection numbers, normalized over the residue field, are unchanged, so $\pi_K$ is minimal and its dual graph agrees with that of $\pi$.
Theorem \ref{Introthm:classification of SFR} now gives strong $F$-regularity of $X_K$.
\end{proof}

\section*{Declaration of AI use}
The manuscript was written entirely by the authors. ChatGPT Astra (OpenAI) was subsequently used to check the mathematical arguments and to assist with language editing. Minor revisions were made based on its suggestions. The authors have carefully reviewed all AI-assisted revisions and take full responsibility for the content of the manuscript.

\section*{Acknowledgements}
We wish to express our gratitude to Kenta Hashizume, Hiroki Matsui, Austyn Simpson, and Shou Yoshikawa
for valuable conversations.
Kawakami was supported by JSPS KAKENHI Grant Numbers JP22KJ1771 and JP24K16897, and by the Inamori Foundation.
Sato was supported by JSPS KAKENHI Grant Number JP24K16900.

\bibliography{references}
\bibliographystyle{alpha}

\bigskip

\end{document}